\documentclass{article}
\usepackage[english]{babel}
\usepackage[utf8]{inputenc}

\usepackage{bm}
\usepackage{amsmath}
\usepackage{graphicx}
\usepackage{color}
\usepackage{multirow}
\usepackage{subcaption}
\usepackage{amsfonts}
\graphicspath{ {./figs/} } 

\date{}
\title{Multiscale model reduction technique for fluid flows with heterogeneous porous inclusions}

\author{
Maria Vasilyeva, S. M. Mallikarjunaiah and D. Palaniappan \\
Department of Mathematics and Statistics\\ Texas A\&M University - Corpus Christi\\
6300 Ocean Drive, Corpus Christi, TX - 78412, USA}

\begin{document}

\maketitle

\begin{abstract}
Numerical treatment of the problem of two-dimensional viscous fluid flow in and around circular porous inclusions is considered.  The mathematical model is described by Navier-Stokes equation in the free flow domain $\Omega_f$ and nonlinear convective Darcy-Brinkman-Forchheimer equations in porous subdomains $\Omega_p$.  
It is well-known that numerical solutions of the  problems in such heterogeneous domains require a very fine computational mesh that resolve inclusions on the grid level.  
The size alteration of the relevant system requires model reduction techniques. 
Here, we present a multiscale model reduction technique based on the Generalized Multiscale Finite Element Method (GMsFEM).  We discuss construction of the multiscale basis functions for the velocity fields based on the solution of the local problems with and without  oversampling strategy. 
Three test cases are considered for a given choice of the three key model parameters, namely, the Reynolds number ($Re$), the Forchheimer coefficient ($C$) and the Darcy number ($Da$). For the test runs, the Reynolds number values are taken to be $Re = 1, 10, 100$ while the Forchheimer coefficient and Darcy number are chosen as $C= 1, 10$ and $Da = 10^{-5}, 10^{-4}, 10^{-3}$, respectively. 
We numerically study the convergence of the method as we increase the number of multiscale basis functions in each domain, and observe good performance of the multiscale method.
\end{abstract}

\section{Introduction}

Mathematical modeling and simulation of fluid flows in the presence of a single or multiple obstacles has been a topic of interest for several decades due to their wide applicability in many practical circumstances across various disciplines. Flow past solid bodies such as cylinders and airfoils have been investigated broadly for a long time (see for instance, \cite{Fornberg1980,Li1995,Sohankar1998,Anderson1998,Yu2009}) by using Navier-Stokes model equations. The presence of nonlinear convective terms poses greater challenges for mathematicians and numerical scientists to solve the PDEs with field boundary and initial conditions. Despite this, several studies have been performed and achieved reasonable theoretical results. The situation is more challenging in the case of fluid flows through and around permeable objects. Indeed, when a flow encounters a porous/permeable object or a collection of objects, a complex flow field develops partially through and partially around the object. In such situations, the prediction of flows (including velocity and pressure fields) passing through and around the object is not straightforward and depends on many factors of the medium in question. Newtonian flows in the presence of a porous cylinder of circular and diamond cross-sections have been analysed recently using numerical methods such as finite volume method (FVM)~\cite{Alazmi2001,Chen2009,rashidi2014numerical,Valipour2014}. The methods used in those studies become more difficult to apply in the presence of numerous inclusions in the flow field. In this study, we investigate two-dimensional viscous incompressible fluid flow in and around multiple circular porous cylinders of arbitrary radii and occupying different positions. We consider Navier-Stokes equations in the exterior flow domain (clear fluid domain) and Darcy-Brinkman-Forchheimer model equations in the porous subdomains and solve the boundary value problem numerically using our proposed multiscale method. 
      
Our primary goal is to devise powerful numerical schemes that can handle the complex nonlinear equations in the exterior and interior porous fluid domains. In order to solve the equations in the porous domains, we consider multiscale model reduction technique for solution of the two dimensional Convective Darcy-Brinkman-Forchheimer (CDBF) equations.  
The nonlinear CDBF equations characterize the motion of incompressible fluid flows in a saturated porous medium and used when the flow velocity is too large for the Darcy’s law to be valid alone \cite{nield2006convection, brinkman1949calculation, mikelic2000homogenization}.   
In \cite{brinkman1949calculation} the set of equations later known as Brinkman equations, intermediate between the Darcy and Stokes equations, is introduced.  
The Brinkman-type law  can be derived using homogenization of the Navier-Stokes equations in a domain containing many tiny solid obstacles \cite{allaire1991homogenization}. 
The numerical study of the two-dimensional fluid flow and forced convection heat transfer around and through a square diamond-shaped porous cylinder is considered in \cite{rashidi2014numerical}. 
Numerical solution is performed using finite volume method. The study pointed out the effects of Reynolds and Darcy numbers on the flow structure and heat transfer characteristics. 
The multiscale model for incompressible fluid flow in porous media with fractures, based on the coupling of Darcy and Brinkman equations is presented in \cite{lesinigo2011multiscale}. The authors provided a finite element scheme for the approximation of the coupled problem, and discuss solution strategies.  The numerical results related to several test cases highlight the potential of this model to reproduce the relevant aspects related to the presence of fractures.
In \cite{Durlofsky_Br} the general methodology of Stokesian dynamics is applied to determine the form of the fundamental solution for flow in porous media. The authors show that the system for dilute porous media behaves as a Brinkman medium.

The problem under consideration involves a global domain with many small porous obstacles (subdomains) that can have various sizes and arbitrary locations. 
The fluid flow in such heterogeneous media have multiscale nature and numerical solution is expensive due to the mesh resolution.  Therefore, viable model reduction methods are necessary in order to improve the computational efficiency and solve the problem on a coarser mesh grid which has much larger length scale compared to the size of perforations. 
Many model reduction techniques, such as numerical homogenization, upscaling and multiscale methods \cite{allaire1992homogenization, allaire1991homogenization, popov2009multiscale, brown2013efficient, mikelic2000homogenization, arbogast2006homogenization, efendiev2013generalized, efendiev2009multiscale, ReducedCon} have been proposed in the literature. 
The upscaling method for solution of the Stokes- Brinkman equation is given in \cite{popov2009multiscale} with applications to naturally fractured karst reservoirs. The Stokes-Brinkman linear model has been used to represent a porous media with a free flow region (fractures, vugs, caves) as a single system of equations. The cell problems that are needed to compute coarse-scale permeability of Representative Element of Volume (REV) are discussed in the cited reference.   
In \cite{araya2017multiscale}, a multiscale hybrid-mixed method is presented to the Stokes and Brinkman equations with highly heterogeneous coefficients. 
A mixed generalized multiscale finite element method for solution of the two dimensional Brinkman equations in the presence of high-contrast permeability fields is discussed in \cite{galvis2015generalized}. The work reported the stability of the mixed multiscale method along with the derivations of a priori error estimates. A variety of two-dimensional numerical examples are also presented to illustrate the effectiveness of the algorithm. 
Another mixed finite-element method in which the Stokes-Brinkman equations are used to compute basis functions for Darcy-flow model on a coarse scale is presented in \cite{gulbransen2010multiscale}. The authors obtained numerical results for strongly heterogeneous sandstone reservoirs, and models of fractured and vugular media.
Further, in \cite{iliev2011variational} a two-scale finite element method for solving Brinkman’s and Darcy’s equation is offered.  The method uses a discontinuous Galerkin finite element method and the concept of subgrid approximation as in \cite{arbogast2006subgrid}.  The proposed algorithms are implemented using the Deal.II finite element library and are tested for a number of model problems.

The fluid flow problem around and through multiple porous cylinders for low and moderate Reynolds numbers and with low and high Forchheimer numbers does not seem to have been addressed adequately in the literature. This is exactly the focus of the present work. As a solution technique, a single-global-continuum-domain approach is assumed which contains the porous cylinders as subdomains. This leads to a single two-dimensional momentum equation, namely, modified Navier-Stokes equations with an additional Darcy and Forchheimer terms. We use a discontinuous Galerkin finite element method to construct fine and coarse grid approximations  \cite{DG_unified, riviere2008discontinuous, girault2005discontinuous, Laz_Tom_Vas_2001}.  We assume that the inclusions have various sizes and arbitrary locations. The macroscopic equations are formulated on a coarse grid with mesh size independent of the size of perforations.  We extend multiscale approach presented in our previous papers \cite{chung2017conservative, alekseev2019multiscale, vasilyeva2021multiscale} to solve the nonlinear convective Darcy-Brinkman–Forchheimer equations in heterogeneous domains.  The multiscale solver for the coarse grid approximation is constructed using Generalized Finite Element Method (GMsFEM) \cite{efendiev2013generalized}. In GMsFEM, we generate a set of multiscale basis functions by construction of the snapshot space and solution of the local spectral problems to reduce the size of the snapshot space. We then solve the equations numerically for three test cases and compute errors by varying the number of multiscale basis functions.

The paper is organized as follows.  The two-dimensional flow past multiple cylinders problem is formulated in section 2. Fine-scale approximation set up along with a discussion on variational and discrete formulations are also provided in the same section.  In Section 3, we present the multiscale method for nonlinear CDBF equations in heterogeneous domains with details. Our numerical results for various test runs including physical descriptions of the numerical solutions and the error calculations are recorded in Section 4. A discussion on typical flow patterns emerging from our numerical computations is also given. The paper ends with some noteworthy points of the present numerical study in conclusion Section 5.

\section{Computational formulation}

\begin{figure}[h!]
\centering
\includegraphics[width=0.6\linewidth]{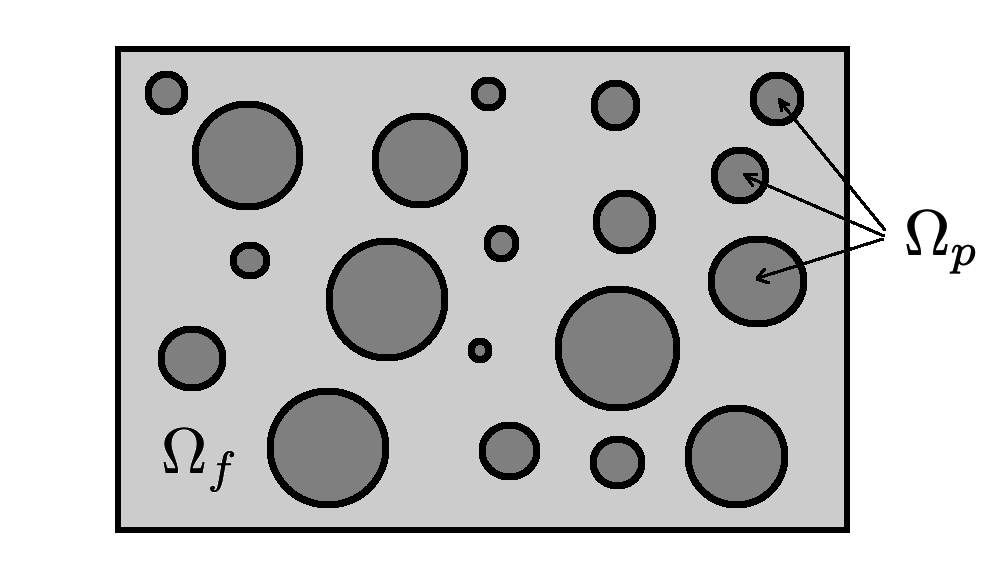}
\caption{Illustration of the computational domain $\Omega = \Omega_f \cup \Omega_p$, where $\Omega_f$ is the Navier-Stokes flow subdomain and  $\Omega_p$ is the subdomain with convective Brinkman-Forchheimer flow.}
\label{sch}
\end{figure}

Consider the problem of a two-dimensional flow of a viscous incompressible fluid through and around multiple permeable cylinders of circular cross-section. A typical computational domain in the problem under consideration is illustrated in Figure \ref{sch}. As in this figure, let $\Omega = \Omega_f \cup \Omega_p$ denote the entire domain consisting of the free fluid region designated by $\Omega_f$ and the flow within porous inclusions represented by $\Omega_p$, respectively. The positions and radii of the circular inclusions are taken to be arbitrary. The assumptions on the governing field equations in the two fluid regions and the interface conditions for the mathematical boundary value problem are as follows.
\begin{itemize}
\item In the free flow subdomain $\Omega_f$, incomprehensible, time-dependent Navier-Stokes equations are assumed.
\item The flow inside the porous subdomain $\Omega_p$ is taken to be that described by convective Darcy-Brinkman-Forchheimer (DBF) model~\cite{popov2009multiscale, popov2009multiphysics,  spiridonov2019mixed, he2021generalized}.
\item At the interface separating domains $\Omega_f$ and $\Omega_p$, we assume stress and velocity continuity boundary conditions.
\end{itemize}

With the above considerations, the governing equations in non-dimensional forms can be provided in the two subdomains. The dimensionless form of the Navier-Stokes equations in the subdomain $\Omega_f$ are   
\begin{equation}
\label{mf}
\begin{split}
\frac{\partial u}{\partial t} + u \cdot \nabla u + \nabla p - \frac{1}{Re} \Delta u = 0, \quad x \in \Omega_f,\\
\nabla \cdot u = 0,  \quad x \in \Omega_f,
\end{split}
\end{equation}
where $u$ represents the two-dimensional velocity field, $p$ the pressure, and 
$Re = \frac{\rho U L}{\mu}$ is the Reynolds number of the flow of the fluid with density $\rho$ and characteristic velocity $U$ in an environment with typical dimension $L$. 
In porous/permeable sub-region $\Omega_p$, we consider the convective Darcy-Brinkman-Forchheimer equations in the form
\begin{equation}
\label{mp}
\begin{split}
\frac{1}{\phi^2} \left( 
\phi \frac{\partial u}{\partial t} 
+ u \cdot \nabla u \right) 
+ \nabla p - \frac{1}{\phi Re} \Delta u 
+ \frac{1}{Re Da} u 
+ \frac{C}{\sqrt{Da}} |u| u
= 0, \quad x \in \Omega_p,\\
\nabla \cdot u = 0,  \quad x \in \Omega_p,
\end{split}
\end{equation}
where 
$Da = \frac{K}{L^2}$ is the Darcy number with permeability $K$, 
$C$ is the Forchheimer coefficient, and 
$\phi$ is the porosity.

Let us define
\[
\xi = 
\left\{\begin{matrix}
1 &  \text{ in } \Omega_f \\
\phi & \text{ in } \Omega_p 
\end{matrix}\right. , \quad 
\chi  = 
\left\{\begin{matrix}
0 &  \text{ in } \Omega_f \\
1 & \text{ in } \Omega_p 
\end{matrix}\right. .
\]
Then equations \eqref{mf} and \eqref{mp} can be written in a combined form as 
\begin{equation}
\label{mm}
\begin{split}
\frac{1}{\xi^2} \left( 
\xi \frac{\partial u}{\partial t} 
+ u \cdot \nabla u \right) 
+ \nabla p 
- \frac{1}{\xi Re} \Delta u 
+ \chi  \left( \frac{1}{Re Da} u 
+ \frac{C}{\sqrt{Da}} |u| u \right)
= 0, \quad x \in \Omega,\\
\nabla \cdot u = 0,  \quad x \in \Omega,
\end{split}
\end{equation}
with  zero initial condition at $t = 0$ and simulate for $t < T_{max}$. 
The interface boundary conditions are written as
\begin{equation}
\label{bc}
u = g_D, \quad x \in \Gamma_u,\quad 
(\nabla u - pI)n = 0, \quad x \in\Gamma_p, 
\end{equation}
where $\Gamma_u \cup \Gamma_p = \partial \Omega$, $n$ is the unit outward normal vector on $\partial \Omega$ and $I$ is the $d\times d$ identity matrix. 

\noindent
\textbf{Variational formulation.}
For the numerical solution of the problem given in \eqref{mm}, we use an implicit scheme for time approximation with linearization from previous time layer.   
To this end, let 
$V = (H^1(\Omega))^2$ and $Q = L^2(\Omega)$, then variational formulation of \eqref{mm} is given by: find $(u^{l+1}, p^{l+1})$ in  $V \times Q$ such that 
\begin{equation}
\label{vf}
\begin{split}
m(u^{l+1}, v) +  a(u^{l+1},v) + d(u^{l+1},v) 
+ b(v, p^{l+1}) &= m(u^l, v), \quad \forall v \in V\\
 b(u^{l+1},q) &= 0, \quad \forall q \in Q
\end{split}
\end{equation}
with
\[
m(u, v) = \frac{1}{\tau} \int_{\Omega} \frac{1}{\xi} \ u \cdot v \ dx, \quad 
a(u, v) = \int_{\Omega} \left(  \frac{1}{\xi Re} \nabla u \colon \nabla v +  \frac{1}{\xi^2} \ (u^l \cdot  \nabla u) \cdot v \right)  dx, 
\]\[
d(u, v) = \chi \int_{\Omega} \left( \frac{1}{Re Da} u \cdot v + \frac{C}{\sqrt{Da}} |u^l| u  \cdot v \right) dx, \quad 
b(v, q) = - \int_{\Omega} q \, \nabla \cdot v \ dx, 
\]
where $l$ is the time layer, $\tau = T_{max}/L$ is the time step and $L$ is the number of time steps, $l = 1,...,L$.

\noindent 
\textbf{Discrete problem.}
For the approximation by space, we use a discontinuous Galerkin method (Interior Penalty Discontinuous Galerkin, IPDG) \cite{DG_unified, riviere2008discontinuous, girault2005discontinuous, Laz_Tom_Vas_2001}.   

\begin{figure}[h!]
\centering
\includegraphics[width=0.55\linewidth]{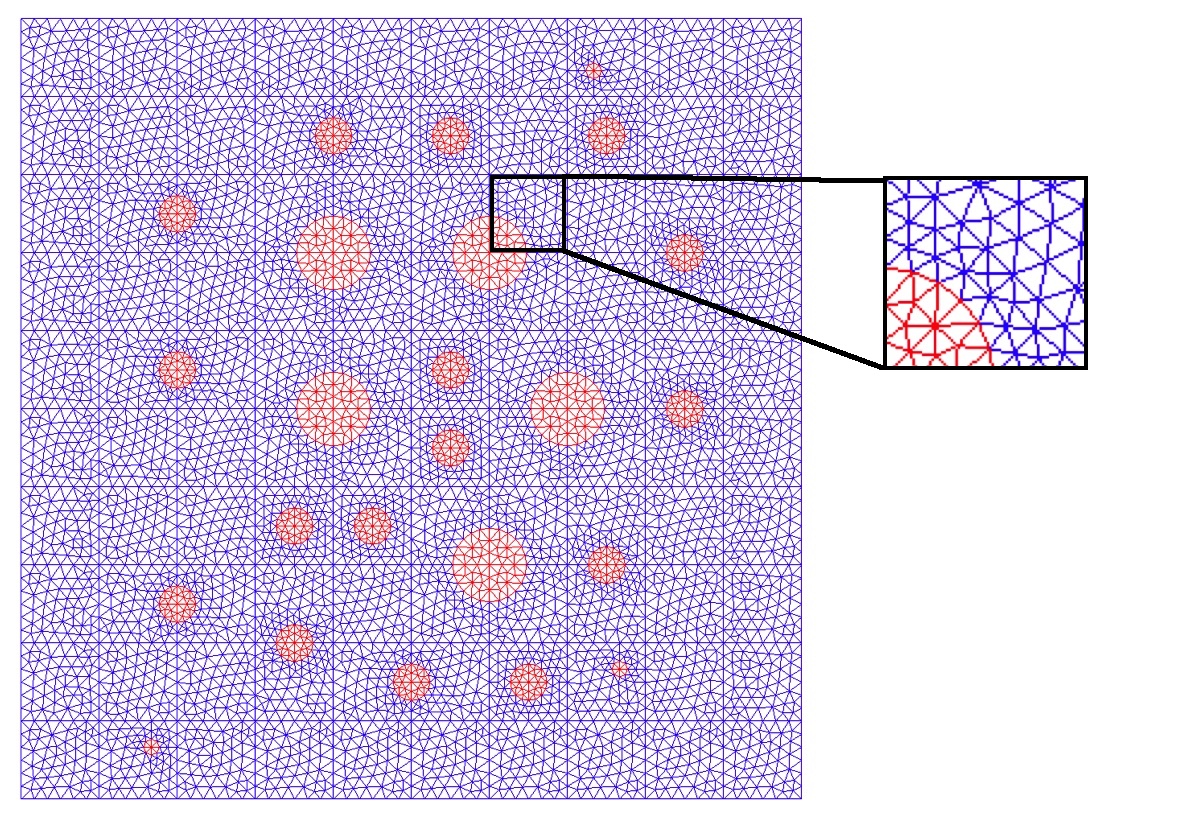}
\caption{Computational domain with circle inclusions and fine grid, $\mathcal{T}^h$.  Blue color: free flow domain. Red color: porous inclusions }
\label{meshf}
\end{figure}

Let $\mathcal{T}^h$ be a fine-grid partition of the domain $\Omega$ that resolve porous inclusions on the grid level with mesh size $h$ (see Figure \ref{meshf}).   
We use the notations $K$ and $E$ to denote a cell and an edge in $\mathcal{T}^h$. 
Let $\mathcal{E}^h$ be the set of edges in $\mathcal{T}^h$ and $\mathcal{E}^h = \mathcal{E}^h_{int} \cup \mathcal{E}^h_{out}$,
where $\mathcal{E}^h_{int}$ is the set of interior edges
and $\mathcal{E}^h_{out}$ is the set of boundary edges.  
For each interior edge $E\in \mathcal{E}^h_{int}$, we define the jump $[u]$ and the average $\{u\}$ of a function $u$ by
\[
[u]_E = u|_{K^{+}} - u|_{K^{-}}, \;\; \{u\}_E = \frac{u|_{K^{+}} + u|_{K^{-}}}{2},
\]
where $K^+$ and $K^-$ are the two  elements sharing the edge $E$,
and the unit normal vector $n$ on $E$ is defined so that $n$ points from $K^+$ to $K^-$.

The fine scale velocity space 
$V_h = \{v_h \in (H^1(\Omega))^2: \,
v_h|_{K} \in (\mathbb{P}_1(K))^2, \, \forall K \in \mathcal{T}^h\}$  contains functions which are piecewise linear in each fine-grid element $K$ and  discontinuous across coarse grid edges.  For the pressure, we use the space of piecewise constant functions $Q_h$.  
We now have following variational formulation for the flow problem \eqref{vf} \cite{riviere2008discontinuous}: find $(u_h^{l+1},p_h^{l+1}) \in V_h \times Q_h$ such that 
\begin{equation}
\label{fau}
\begin{split}
m(u_h^{l+1}, v_h)+a_{DG}(u_h^{l+1}, v_h)+d(u_h^{l+1}, v_h)+b_{DG}(p_h^{l+1},v_h)
&=m(u_h^l, v_h) + f(v_h), \quad \forall v_h \in V_h\\
b_{DG}(u_h^{l+1}, q_h)&= l(q_h), \quad \forall q_h \in Q_h,    
\end{split}
\end{equation}
where
\[
\begin{split}
m(u,v) &=   \frac{1}{\tau} \sum_{K \in \mathcal{T}^h} \int_K  \frac{1}{\xi} \ u \cdot v \ dx,
\\
a_{DG}(u,v) &=   \sum_{K \in \mathcal{T}^h} \int_K  \left(  \frac{1}{\xi Re} \nabla u \colon \nabla v  +  \frac{1}{\xi^2} \ (u^l \cdot  \nabla u) \cdot v \right)  dx\\
  &- \sum_{E \in \mathcal{E}^h_{int}} 
\int_E  \left(  
    \left\{ \frac{1}{\xi Re} \nabla u \right\} n \cdot [v]
+  \left\{ \frac{1}{\xi Re} \nabla v \right\} n \cdot [u] 
- \frac{\gamma}{\{\xi Re \ h\}} [u] \cdot [v]
\right) \, ds\\
&- \sum_{E \in \mathcal{E}^h_{out, D}} 
\int_E  \frac{1}{\xi Re} \left(  
(\nabla u \ n ) \cdot  v + (\nabla v  \ n ) \cdot u - 2 \frac{\gamma}{h} u \cdot v
\right) \, ds,  
\end{split}
\]\[
\begin{split}
d(u, v) &= \chi \sum_{K \in \mathcal{T}^h} \int_K \left( \frac{1}{Re Da}u \cdot v + \frac{C}{\sqrt{Da}} |u^l| u \cdot v \right) dx,
\\
b_{DG}(v, q) &= - \sum_{K \in \mathcal{T}^h} \int_K q \ \nabla \cdot v \ dx
 + \sum_{E \in \mathcal{E}^h_{int}} \int_E \{q\} [v] \cdot n \ ds 
 + \sum_{E \in \mathcal{E}^h_{out, D}} \int_E q \ v \cdot n  \ ds, 
\end{split}
\]\[
\begin{split}
f(v)  &
=  \sum_{E \in \mathcal{E}^h_{out, D}} 
\int_E  \frac{1}{\xi Re} \left(  
 (\nabla g_D \ n ) \cdot v +  (\nabla v  \ n ) \cdot g_D  - 2 \frac{\gamma}{h } g_D \cdot v
\right) \, ds,
\\
l(q)  &= \sum_{E \in \mathcal{E}^h_{out, D}} \int_E q \ g_D \cdot n  \ ds.
\end{split}
\]
Here $\gamma$ is the penalty parameter, $n$ is the unit normal to the edge $E$ and $\mathcal{E}^h_{out, D}$ is the set of boundary edges related to the boundary $\Gamma_D$.
One can see that $(u_h, p_h)$ will converge to the exact solution $(u, p)$ in the energy norm as the fine mesh size $h \rightarrow 0$. 
Moreover the fine mesh is constructed to resolve interface between two subdomains $\Omega_f$ and $\Omega_p$ (see Figure \ref{meshf}). 

We can write the above discrete systems in the matrix form as follows. 
\begin{equation}
\begin{split}
\begin{pmatrix}
M_h + A_h + D_h & B_h^T \\
B_h & 0
\end{pmatrix}
\begin{pmatrix}
u^{l+1}_h \\
p^{l+1}_h
\end{pmatrix} = 
\begin{pmatrix}
M_h u^l_h + F_h^u \\
F_h^p.
\end{pmatrix}.
\end{split}
\end{equation} 
where
\[
M_h = [m_{ij} = m(\psi_i, \psi_j)],  \quad 
A_h = [a_{ij} = a_{DG}(\psi_i, \psi_j)],  \quad 
D_h = [d_{ij} = d(\psi_i, \psi_j)],  
\]\[
B_h = [b_{ij} = b_{DG}(\phi_i, \psi_j)],  \quad 
F_h^u = [f_{j} = f(\psi_j)],   \quad
F_h^p = [f_{j} = l(\phi_j)],   
\]
with a piecewise linear function  $\psi_i  \in (\mathbb{P}_1(K))^2$  in each fine-grid element $K \in \mathcal{T}^h$  is and  discontinuous across coarse grid edges.  For the pressure,  $\phi_i$  is the  piecewise constant functions on mesh $\mathcal{T}^h$. The size of the discrete system is $DOF_h = 3 \cdot d \cdot N_h + N_h$, where $d$ is the dimension ($d = 2$) and $N_h$ is the number of fine grid cells. 

In the following section, we provide our multiscale method used for the size reduction of the system described above. 
In the adopted multiscale method, we solve problems in local domains with various boundary conditions to form a snapshot space and use a spectral problem in the snapshot space to perform the required dimension reduction.

\section{Multiscale method}

We begin by describing construction of the coarse grid approximation using Generalized Finite Element Method (GMsFEM) \cite{alekseev2019multiscale,  ceh2016adaptive, AdaptiveGMsDGM}.  

\begin{figure}[h!]
\centering
\includegraphics[width=0.34\linewidth]{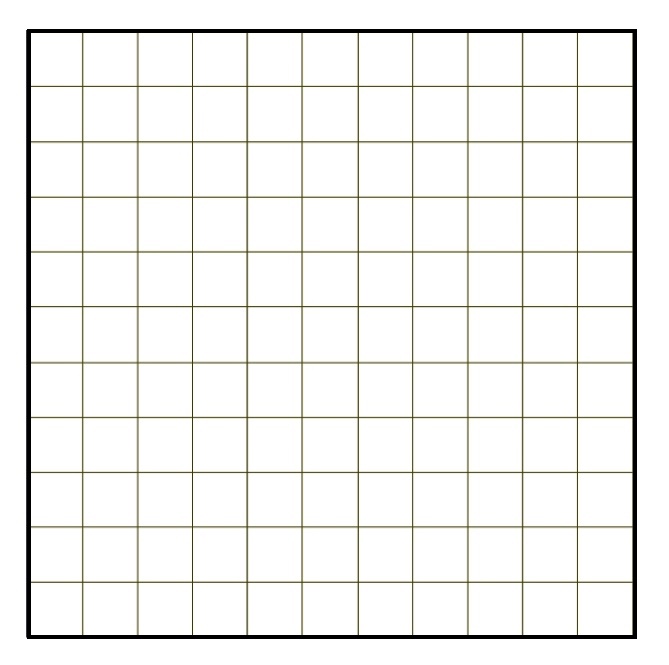}
\includegraphics[width=0.5\linewidth]{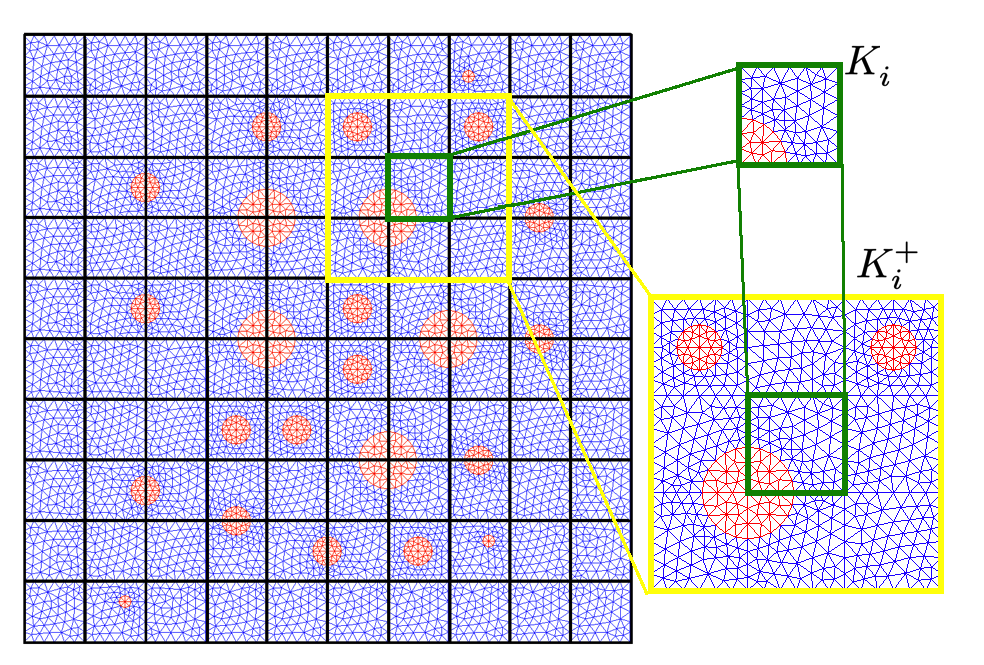}
\caption{Illustration of coarse grid and local domains $K_i$ (coarse cell) and oversampled local domain $K_i^+$. 
Left: $10 \times 10$ coarse grid.  
Right: Local domain with fine grid resolution that conformed with inclusions and coarse grid facets}
\label{meshc}
\end{figure}

Let $\mathcal{T}^H$ be a coarse grid  of  domain $\Omega$ with mesh size $H$  and $\mathcal{E}^H$ be the set of all facets of the coarse grid in $\mathcal{T}^H$ (see Figure \ref{meshc}).   
For the sake of simplicity, in this work, we consider structured coarse grid with quadratic cells.  In general coarse grid cell can have any shape (unstructured coarse grid) and can be a mesh partitioning \cite{chung2016multiscaledg}.  In the present situation, we consider $10 \times 10$ coarse grid. 
We define $V^{H}$ as the multiscale velocity space, which contains a set of basis functions supported in each coarse block $K$. 
For the pressure approximation, we use the piecewise constant function space $Q^{H}$ over the coarse cells.  

\noindent
We construct a multiscale space for the velocity field
\[
V^{H} = \text{span} \{ \psi_i \}_{i=1}^{N_u}, 
\]
where $N_u =\text{dim}(V^{H})$ is the total number of basis functions.
\noindent
For the pressure, we use the space of piecewise constant functions over the coarse cell
\[
Q^{H} = \{ r \in L^2(\Omega) : \,  r |_K \in P^0(K), \, \forall K \in \mathcal{T}_H \}
\]
where $N_p = \text{dim}(Q^{H})$ and equal to the number of coarse grid cells.  
\noindent
For the coarse grid approximation, we use a Discontinuous Galerkin (DG) approach and have following variational formulation for each time step $l+1$: find $(u^{l+1}_H, p^{l+1}_H) \in V^{H} \times Q^{H}$ such that \cite{MsDG, AdaptiveGMsDGM, chung2016multiscaledg, chung2017conservative}
\begin{equation}
\label{v-ms}
\begin{split}
m(u^{l+1}_H, v_H) +  a(u^{l+1}_H,v_H) + d(u^{l+1}_H,v_H) 
+ b(v_H, p^{l+1}_H) &= m(u^l_H, v_H), \quad \forall v_H \in V^H\\
 b(u^{l+1}_H,q_H) &= 0, \quad \forall q_H \in Q^H.
\end{split}
\end{equation}

The multiscale space for the velocity is constructed by solution of the local spectral problem. We first construct snapshot space for each local domain by the solution of local problem with different boundary conditions. Then, by the solution of the local eigenvalue problem on snapshot space, we select dominant modes of the velocity filed as multiscale basis functions.  

\textbf{Snapshot space.} 
We construct local snapshots in each coarse cell 
$\phi_{i,k} \in K_i$, where  $i = 1, \cdots, N$ and $N$ is the number of coarse cells in $\mathcal{T}^H$. 
The local snapshot space is generated by the solution of the following local problem: 
find  $(\phi_{i,k},  \eta_{i,k}) \in V^h(K_i) \times Q^h(K_i)$  such that
\begin{equation}
\label{eq:snapu}
\begin{split}
 a_{DG}(\phi_{i,k}, v) + d(\phi_{i,k}, v) + b_{DG}(\eta_{i,k}, v) &=  f(v), \quad \forall  v \in V^h(K_i)\\
 b_{DG}(\phi_{i,k}, q) &= l(q), \quad \forall q \in Q^h(K_i). 
\end{split}
\end{equation}
where
\[
\begin{split}
a_{DG}(u,v) &=   \sum_{K \in \mathcal{T}^h(K_i)} \int_K  \left(  \frac{1}{\xi Re} \nabla u \colon \nabla v  +  \frac{1}{\xi^2} \ (\tilde{u} \cdot  \nabla u) \cdot v \right)  dx\\
  &- \sum_{E \in \mathcal{E}^h_{int}(K_i)} 
\int_E  \left(  
     \left\{\frac{1}{\xi Re} \nabla u \right\} n \cdot [v]
+  \left\{ \frac{1}{\xi Re} \nabla v \right\} n \cdot [u] 
- \frac{\gamma}{\{\xi Re \ h\}} [u] \cdot [v]
\right) \, ds\\
&- \sum_{E \in \mathcal{E}^h_{out}(K_i)} 
\int_E \frac{1}{\xi Re} \left(  
(\nabla u \ n ) \cdot v +(\nabla v  \ n ) \cdot u - \frac{ 2 \gamma}{h} u \cdot v
\right) \, ds,  
\end{split}
\]\[
\begin{split}
d(u, v) &= \chi \sum_{K \in \mathcal{T}^h(K_i)} \int_K \left( \frac{1}{Re Da} u \cdot v + \frac{C}{\sqrt{Da}} |\tilde{u}| u \cdot v \right) dx,
\\
b_{DG}(v, q) &= - \sum_{K \in \mathcal{T}^h(K_i)} \int_K q \ \nabla \cdot v \ dx
 + \sum_{E \in \mathcal{E}^h_{int}(K_i)} \int_E \{q\} \ [v] \cdot n \ ds 
 + \sum_{E \in \mathcal{E}^h_{out}(K_i)} \int_E  q \ v \cdot n  \ ds, 
\end{split}
\]\[
\begin{split}
f(v)  &
=  \sum_{E \in \mathcal{E}^h_{out}(K_i)} 
\int_E  \frac{1}{\xi Re}  \left(  
(\nabla g \ n ) \cdot v +   (\nabla v  \ n ) \cdot g - \frac{ 2 \gamma}{ h } g \cdot v
\right) \, ds,
\\
l(q)  &= \sum_{E \in \mathcal{E}^h_{out}(K_i)} \int_E  q \ g \cdot n  \ ds 
+ \sum_{K \in \mathcal{T}^h(K_i)} \int_K c \ q \ dx,
\end{split}
\]
where $g =\delta_i^k$,  $\delta_i^k$ is the discrete delta function defined on $\partial K_i$,  $k = 1, \cdots, J_i$ and $J_i$ is the number of fine grid facets on the boundary of $K_i$. 
Here $c$  is chosen via the compatibility condition, $c = \frac{1}{|K_i|}\int_{\partial K_i} \delta_i^l \cdot n \, ds$.

To handle the nonlinear problem, we use the global solution of linear problem with the given boundary conditions at final time to calculate $\tilde{u}$.  In general,  we can generate a set of a global functions $\tilde{u}_i$ with different boundary conditions  to handle the general case 
\cite{efendiev2014generalized, brown2016generalized, chung2018multiscale}. 
We note that the problem \eqref{eq:snapu} is linear for a given global solution $\tilde{u}$. 
We form a local snapshot space in $K_i$ using all the local solutions 
\[
V^{i,\text{snap}} = \{ \phi_{i,k}: 1 \leq k \leq J_i \}
\] 
and define projection matrix to the snapshot space 
\[
R_{i, \text{snap}} = \left[ \phi_{i,1} \ldots, \phi_{i,J_i} \right]^T.
\]

\textbf{Oversampling strategy in snapshot space construction.} 
In order to reduce the boundary effects in snapshot space construction and  improve the accuracy of  multiscale methods, we apply the oversampling strategy \cite{Efendiev_oversampling13, chung2017conservative, randomized2014}.  

Let $K_i^{+}$ be an enlarged domain of $K_i$ constructed by adding one coarse block around target coarse cell (see  Figure \ref{meshc}). 
To construct snapshot space, we solve the local problem \eqref{eq:snapu} and find  $\phi_{i,k}^{+}$ in an oversampled domain $K_i^{+}$ with similar Dirichlet boundary conditions for the velocity field, $\phi_{i,k}^{+}=\delta_i^k$ on $\partial K_i^{+}$, where $k= 1, \cdots, J_i^{+}$, where $J_i^{+}$ is the number of fine edges on the boundary of $K_i^{+}$.  
Note that the velocity solutions (snapshots) of these local problems are supported in the larger domain $K_i^{+}$. 
To form a local snapshot space in $K_i$, we restrict the functions $\phi_{i,k}^{+}$ on $K_i$ and generate the snapshot basis
\[
V^{i,\text{snap}} = \{ \phi_{i,k}: 1 \leq k \leq J^+_i \},
\] 
where $\phi_{i,k} = \phi_{i,k}^{+}|_{K_i}$.  
We note that the snapshot space contains  extensive number of basis functions and therefore need to employ space reduction technique to form a subspace which can then approximate the snapshot space accurately and consequently improve computational efficiency.

\textbf{Multiscale basis functions for velocity.} 
The size reduction of the snapshot space is achieved by solving the local spectral problem in local domain $K_i$. 
From the following generalized eigenvalue problem in the snapshot space, we find  $(\lambda_k,  \psi^{\text{snap}}_{i,k})$
\begin{equation}
\label{eq:off-eq}
A^{i,\text{snap}} \psi^{\text{snap}}_{i} =  
\lambda  S^{i,\text{snap}} \psi^{\text{snap}}_{i}, 
\end{equation}
where 
\[
A^{i,\text{snap}} = R_{i, \text{snap}} A^i R_{i, \text{snap}}^T, \quad 
S^{i,\text{snap}} = R_{i, \text{snap}} S^i R_{i, \text{snap}}^T. 
\]
Here $A^i$ and $S^i$ are the matrix representation of the bilinear form $a_{DG}(u, v)$ and $s(u, v)$
\[
\begin{split}
a_{DG}(u,v) &=  \sum_{K \in \mathcal{T}^h(K_i)} \int_K   \frac{1}{\xi Re} \nabla u \colon \nabla v \ dx\\
  &- \sum_{E \in \mathcal{E}^h_{int}(K_i)} 
\int_E  \left(  
    \left\{ \frac{1}{\xi Re} \nabla u \right\} \cdot  [v]
+  \left\{ \frac{1}{\xi Re} \nabla v \right\} \cdot [u] 
- \frac{\gamma}{\{\xi Re \ h\}} [u] \cdot [v]
\right) \ ds,\\
s(u,v) &= \sum_{E \in \mathcal{E}^h_{int}(K_i)} 
\int_E u \cdot v \, ds.
\end{split}
\]
Note that the integral in $s(u,v)$ is defined on the boundary of the coarse cell. 
We arrange the eigenvalues in increasing order 
\[
\lambda_1 \leq  \lambda_2 \leq \cdots \leq \lambda_{J_i} 
\]
and choose the first  $M$ eigenvectors corresponding to the smallest eigenvalues as multiscale basis functions for the velocity field. This yields
\[
V^{H} 
= \text{span} \{ \psi_{i,k}: \, 1 \leq i \leq N, \,  1 \leq k \leq M \},  
\]
where 
$\psi_{i,k} = R_{i, \text{snap}} \psi^{\text{snap}}_{i,k}$ and $N$ is the number of coarse cell.  

\textbf{Coarse scale system.} 
We construct the coarse grid system using a global projection approach.  We form the projection  matrices using  the computed multiscale basis functions
\[
R_u = \left[ 
\psi_{1,1} , \ldots, \psi_{1,M}, \ldots, 
\psi_{N,1}, \ldots, \psi_{N, M} \right]^T,
\quad 
R_p = \left[ \eta_1 , \ldots, \eta_{N} \right]^T, 
\]
where $N$ is the number of coarse grid cells and $M$ is the number of multiscale basis functions for the velocity field. 
In general we can apply an adaptive approach and use a different number of the basis functions in each local domain $K_i$, i.e. $M_i$. 
Note that we use the space of piecewise constant functions for pressure over the coarse grid $K_i$, and set $\eta_i(x)$ equal to 1 if $x \in K_i$ and zero otherwise.

Using projection matrices for the velocity and pressure fields, we obtain the following discrete system in matrix form:
\begin{equation}
\label{eq:coarsesystem}
\begin{split}
\begin{pmatrix}
M_H + A_H + D_H & B_H^T \\
B_H & 0
\end{pmatrix}
\begin{pmatrix}
u^{l+1}_H \\
p^{l+1}_H
\end{pmatrix} = 
\begin{pmatrix}
M_H u^l_H + F_H^u \\
F_H^p.
\end{pmatrix}.
\end{split}
\end{equation} 
where 
\[
M_H = R_u M_h R_u^T, \quad 
A_H = R_u A_h R_u^T, \quad 
B_H = R_u B_h R^T_p, \quad 
F_H^u = R_u F_h^u, \quad
F_H^u = R_p F_h^p.
\] 
After solution of the coarse-scale system, we can reconstruct velocity on a fine grid  
\[
u_{ms}= R_u^T u_H.
\]
The size of the resulting discrete system \eqref{eq:coarsesystem} is $DOF_H = N_u + N_p$, where $N_p = N$ and  $N_u = M \cdot N$, where once again, $N$ and $M$ are respectively the number of coarse grid cells and multiscale basis functions for the velocity field.

\section{Numerical results}

\begin{figure}[h!]
\centering
\includegraphics[width=0.4\linewidth]{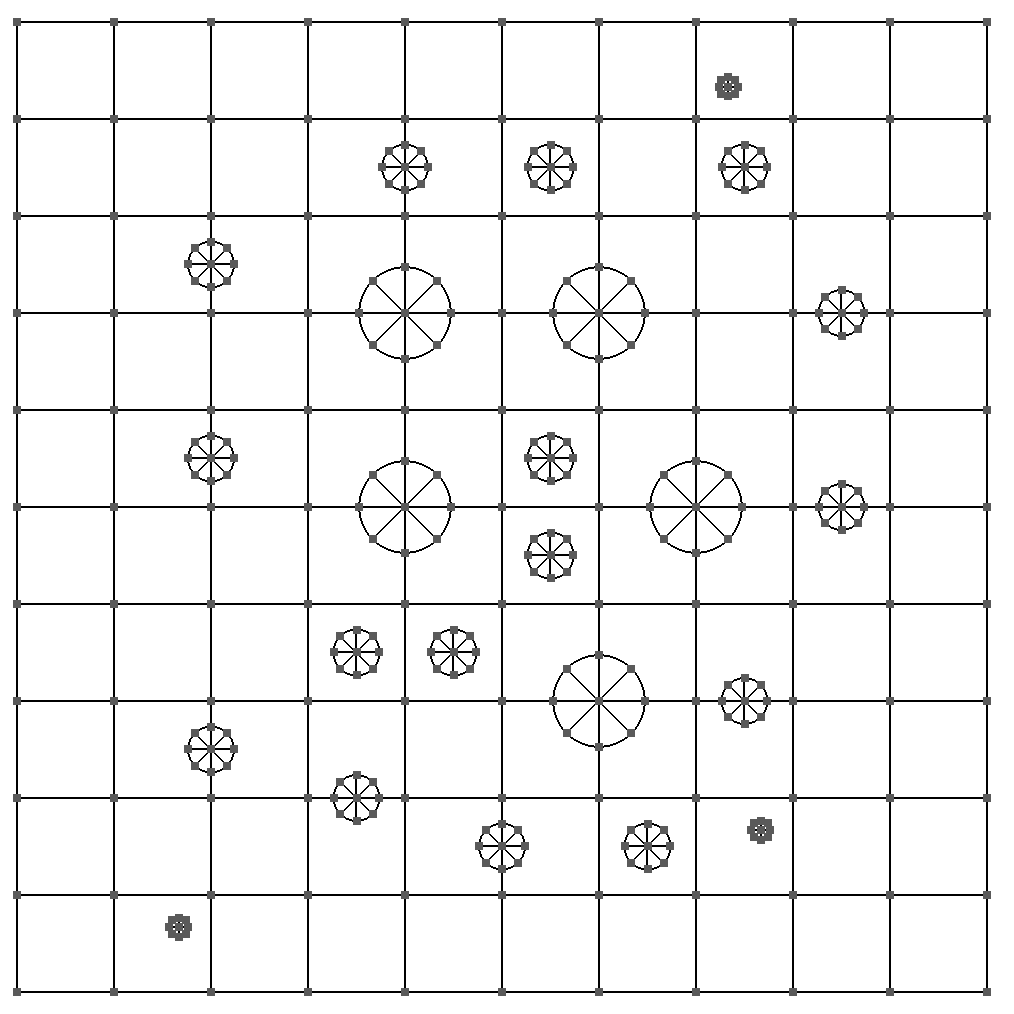}\ \ \ \ \ \ 
\includegraphics[width=0.4\linewidth]{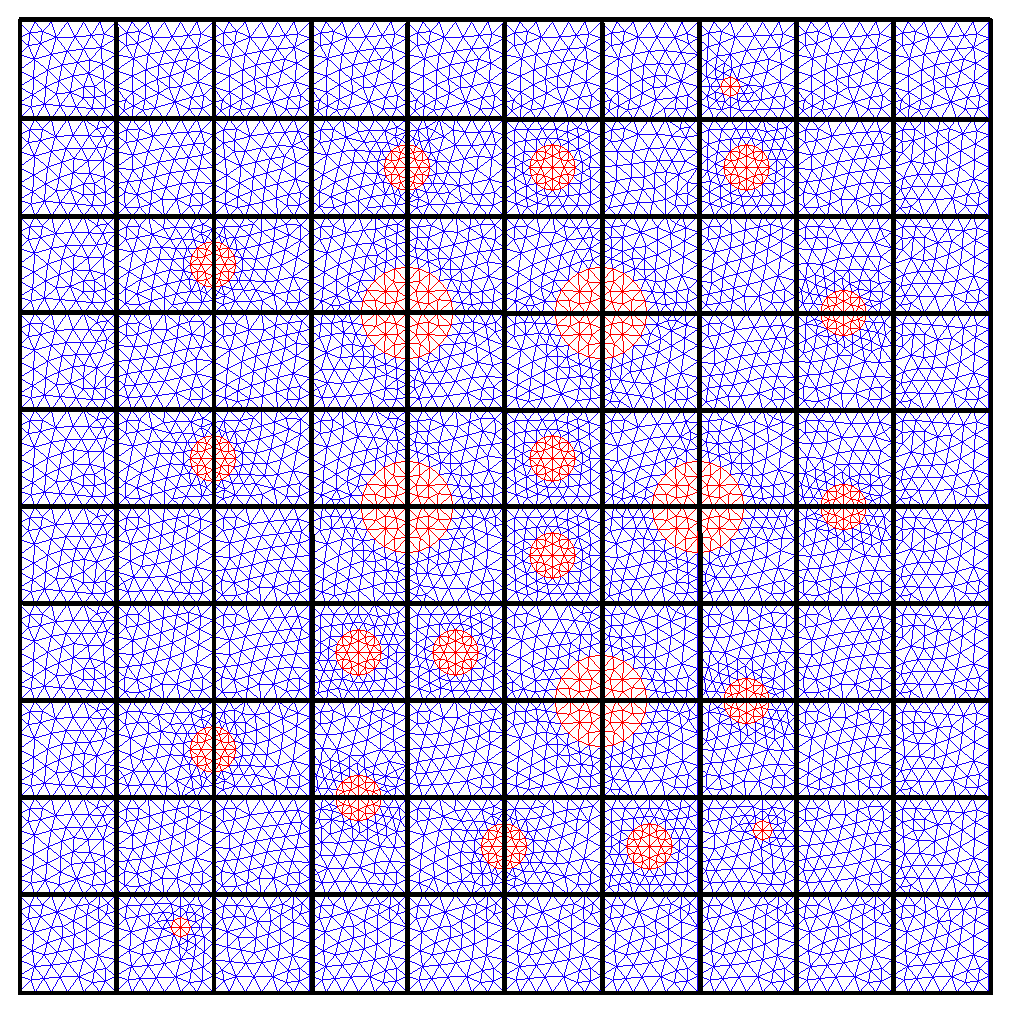}
\caption{Computational domain with circle inclusions (left)  and fine grid with coarse mesh (right)}
\label{mesh}
\end{figure}

We now turn our focus on the numerical solutions and simulations of the nonlinear fluid flow problem in and around multiple circular porous inclusions formulated in section 2.  
The numerical calculations have been performed in the heterogeneous computational domain $\Omega = [-1, 1]^2$ using the multiscale method demonstrated in the previous section. The governing equations \eqref{mm} with boundary conditions \eqref{bc} are solved over two-dimensional fine and coarse grid systems. The fine grid contains  6373 vertices and 12504 triangular cells while the coarse grid size is taken to be $10 \times 10$ in size with 121 vertices and 100 cells. In Figure \ref{mesh}, we have depicted the computational domain and fine grid constructed using Gmsh sofware \cite{geuzaine2009gmsh}.  The fine grid has been built in such a way that it resolves the interface between domain $\Omega_f$ and $\Omega_p$ on a grid level. Moreover, the fine grid is conforming with coarse grid edges, and the coarse grid is uniform with quadratic cells. On the right figure, we have presented fine grid in blue, coarse grid in black, and the circular inclusions in red colors, respectively. Note that we have taken 24 circular inclusions with different radii at relatively random locations.

It is evident from the discussion in Section 2 that there are three key parameters involved in our model problem namely, the Reynolds number $Re$, the Darcy number $Da$, and the Forchheimer coefficient $C$, in addition to the porosity $\phi$. These constants constitute a four-set parameter space in general, however, we restrict our investigation to some fixed set of values to proceed with the numerical computations based on the multiscale method. Specifically, the performance of the multiscale method has been tested for the following three test cases with the given choices of the parameters:
\begin{itemize}
\item \textit{Test 1}: $Re = 1$,  $C = 1$ and $T_{max}=0.01$.
\item \textit{Test 2}: $Re = 10$, $C = 10$ and $T_{max}=0.1$.
\item \textit{Test 3}: $Re = 100$, $C = 1$ and $T_{max}=1.0$.
\end{itemize}
In all our numerical test runs the porosity parameter is set to $\phi=0.3$. Three distinct values for the Darcy number, viz., $Da = 10^{-5}, 10^{-4}$ and $10^{-3}$, are chosen to illustrate the behavior of the velocity and pressure fields for flows around and inside the circular inclusions. For illustration purposes, we have considered 24 circular inclusions. On the left boundary  $\Gamma_u$, we set velocity $u = g_D$ with $g_D = (1,0)$. At the top and bottom boundaries, we set zero velocity and on the right boundary we take zero pressure. We perform numerical simulations with 50 time steps for both fine grid and coarse-grid solutions.

\begin{figure}[h!]
\centering
\begin{subfigure}{0.32\textwidth}
\centering
$Da = 10^{-5}$
\end{subfigure}
\begin{subfigure}{0.32\textwidth}
\centering
$Da = 10^{-4}$
\end{subfigure}
\begin{subfigure}{0.32\textwidth}
\centering
$Da = 10^{-3}$
\end{subfigure}\\
\begin{subfigure}{0.32\textwidth}
\centering
\includegraphics[width=0.78\linewidth]{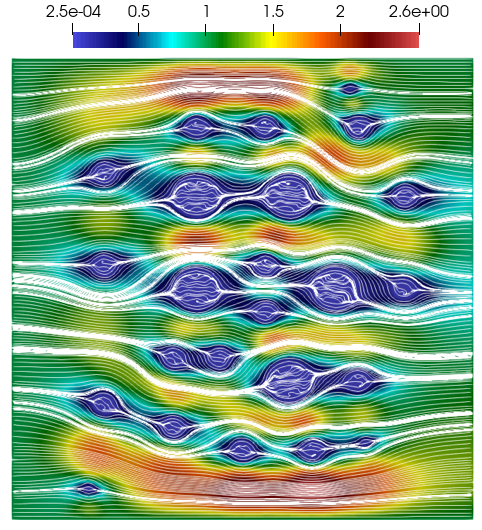}\\
\includegraphics[width=0.78\linewidth]{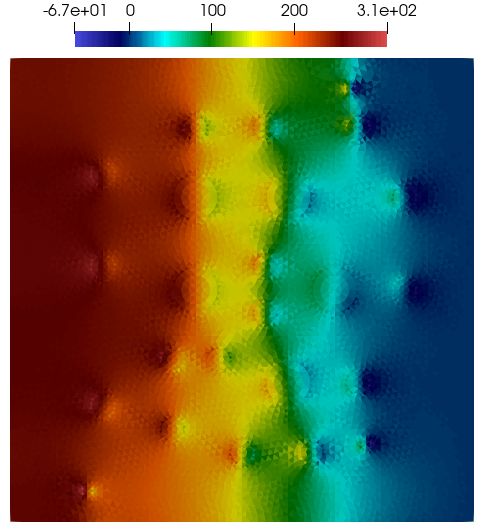}
\end{subfigure}
\begin{subfigure}{0.32\textwidth}
\centering
\includegraphics[width=0.78\linewidth]{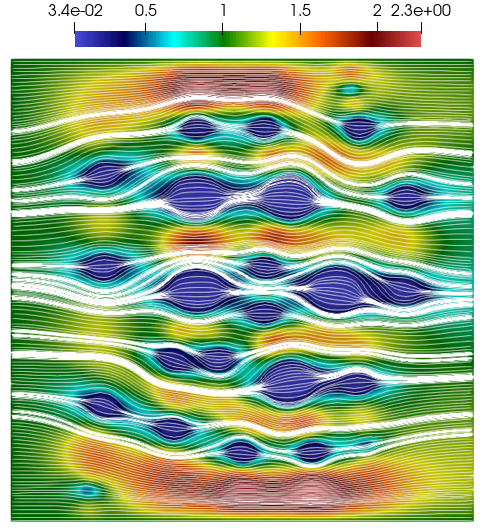}\\
\includegraphics[width=0.78\linewidth]{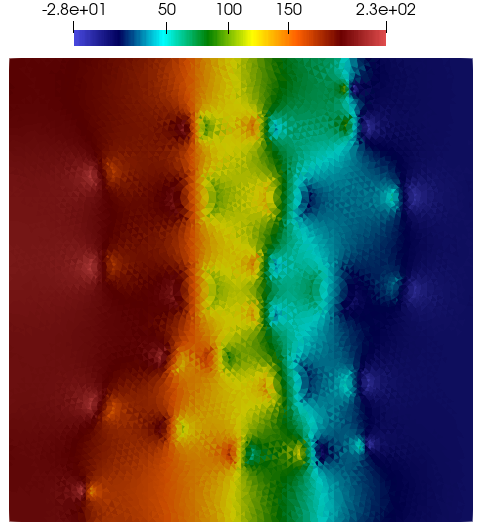}
\end{subfigure}
\begin{subfigure}{0.32\textwidth}
\centering
\includegraphics[width=0.78\linewidth]{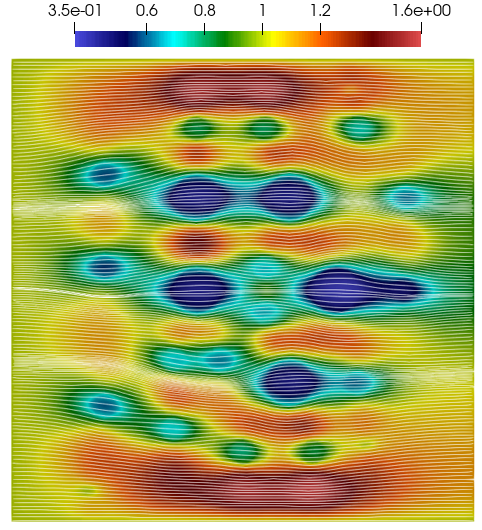}\\
\includegraphics[width=0.78\linewidth]{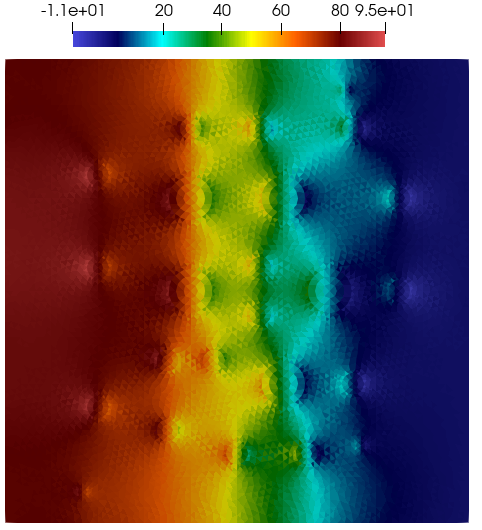}
\end{subfigure}
\caption{\textit{Test 1} ($Re = 1$, $C = 1$).  Fine grid solution at the final time for $Da = 10^{-5}, 10^{-4}, 10^{-3}$ (from left to right).  
First row: velocity magnitude with streamlines. 
Second row: pressure}
\label{up-1}
\end{figure}

\begin{figure}[h!]
\centering
\begin{subfigure}{0.32\textwidth}
\centering
$Da = 10^{-5}$
\end{subfigure}
\begin{subfigure}{0.32\textwidth}
\centering
$Da = 10^{-4}$
\end{subfigure}
\begin{subfigure}{0.32\textwidth}
\centering
$Da = 10^{-3}$
\end{subfigure}\\
\begin{subfigure}{0.32\textwidth}
\centering
\includegraphics[width=0.78\linewidth]{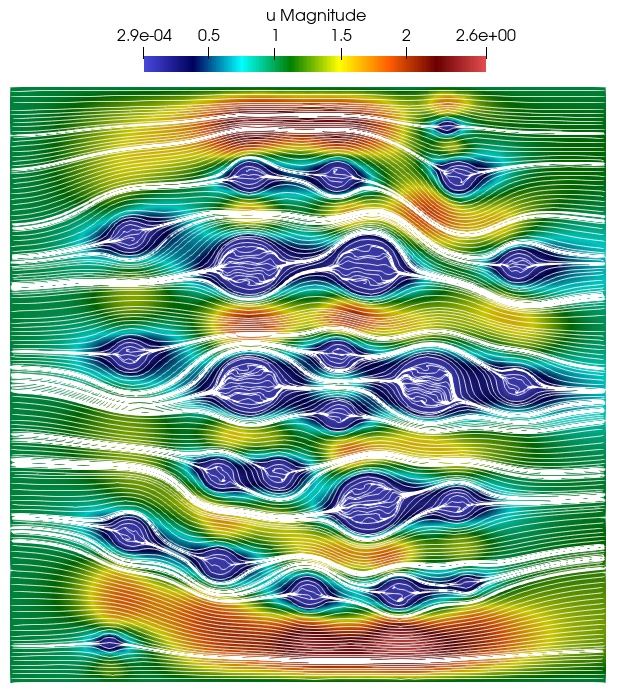}\\
\includegraphics[width=0.78\linewidth]{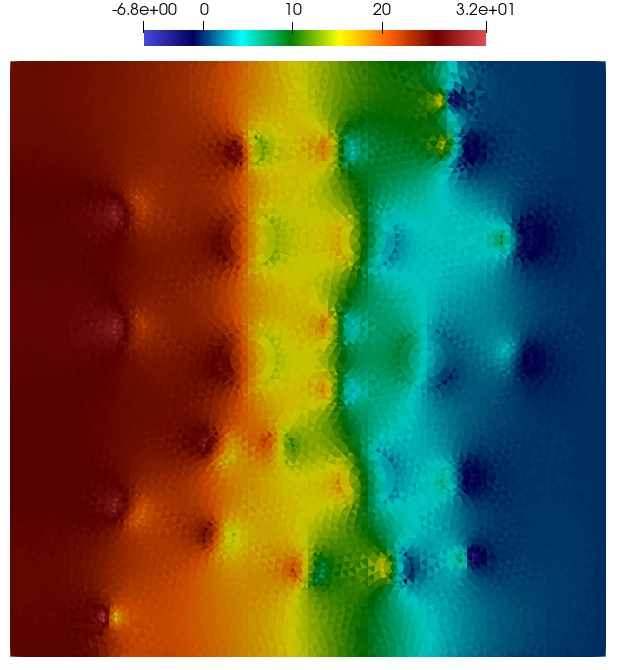}
\end{subfigure}
\begin{subfigure}{0.32\textwidth}
\centering
\includegraphics[width=0.78\linewidth]{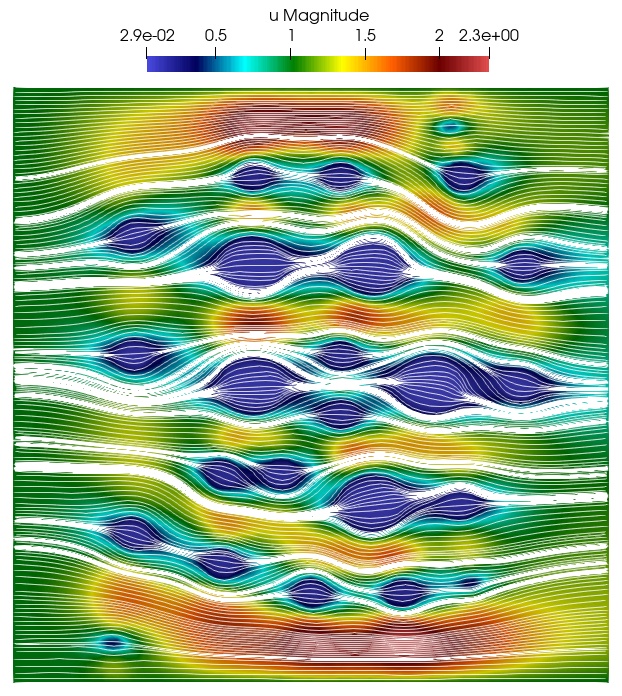}\\
\includegraphics[width=0.78\linewidth]{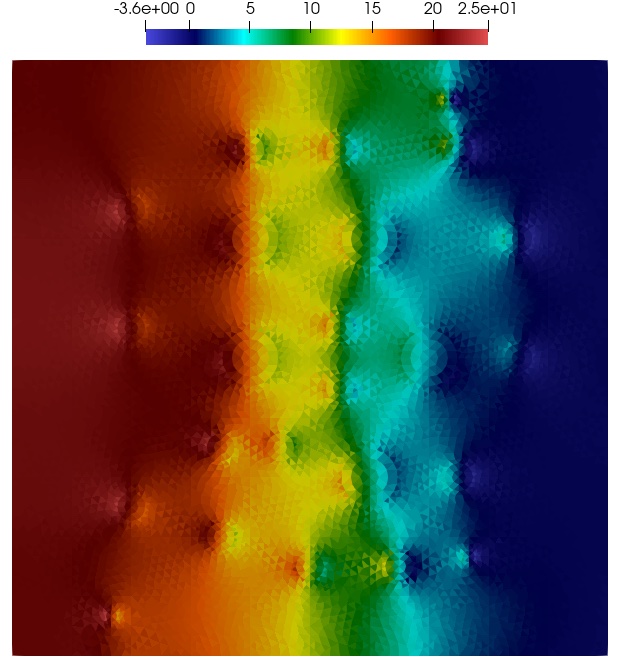}
\end{subfigure}
\begin{subfigure}{0.32\textwidth}
\centering
\includegraphics[width=0.78\linewidth]{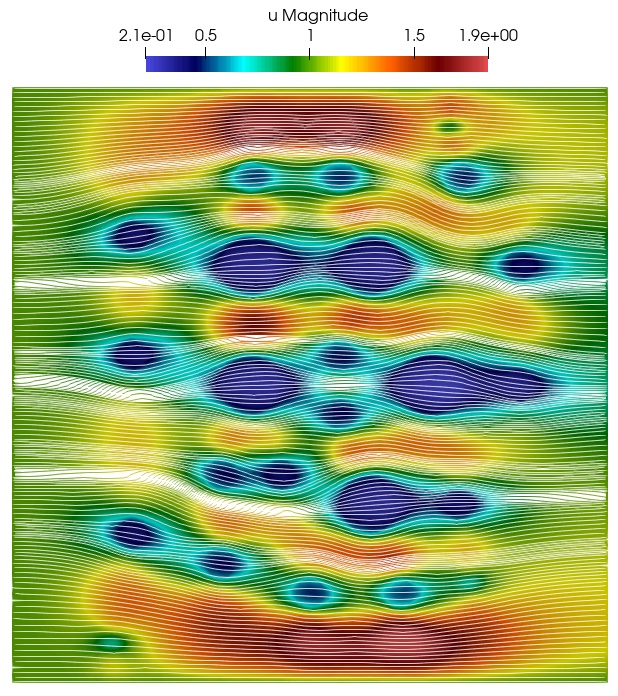}\\
\includegraphics[width=0.78\linewidth]{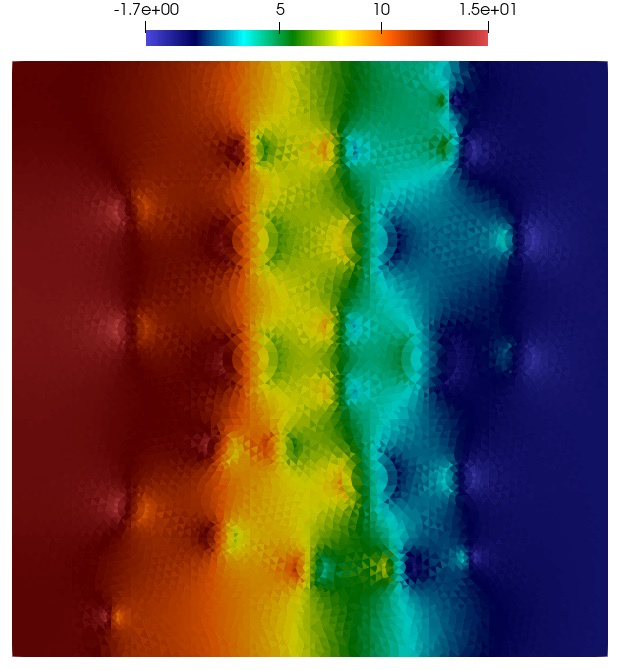}
\end{subfigure}
\caption{\textit{Test 2} ($Re = 10$, $C = 10$).  Fine grid solution at final time for $Da = 10^{-5}, 10^{-4}, 10^{-3}$ (from left to right).  
First row: velocity magnitude with streamlines. 
Second row: pressure}
\label{up-2}
\end{figure}

\begin{figure}[h!]
\centering
\begin{subfigure}{0.32\textwidth}
\centering
$Da = 10^{-5}$
\end{subfigure}
\begin{subfigure}{0.32\textwidth}
\centering
$Da = 10^{-4}$
\end{subfigure}
\begin{subfigure}{0.32\textwidth}
\centering
$Da = 10^{-3}$
\end{subfigure}\\
\begin{subfigure}{0.32\textwidth}
\centering
\includegraphics[width=0.78\linewidth]{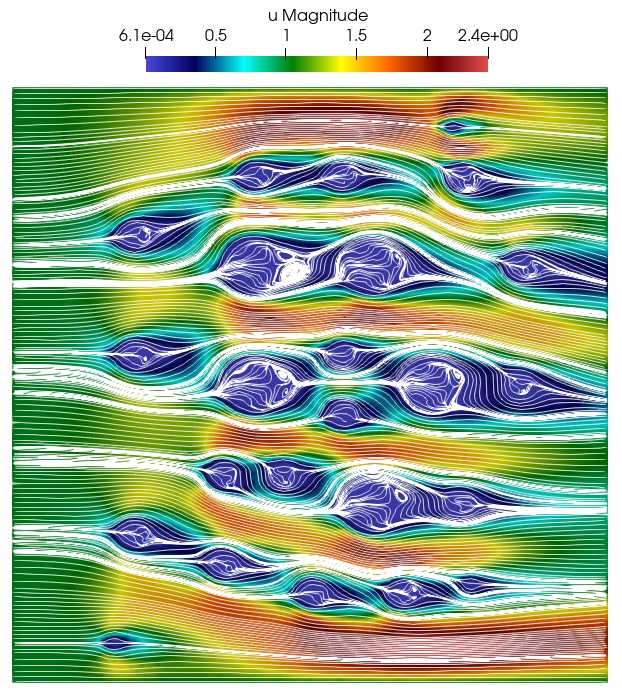}\\
\includegraphics[width=0.78\linewidth]{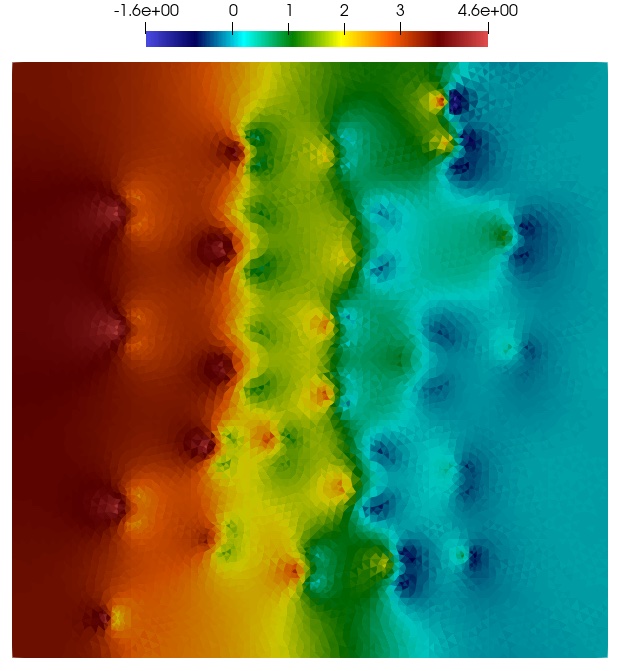}
\end{subfigure}
\begin{subfigure}{0.32\textwidth}
\centering
\includegraphics[width=0.78\linewidth]{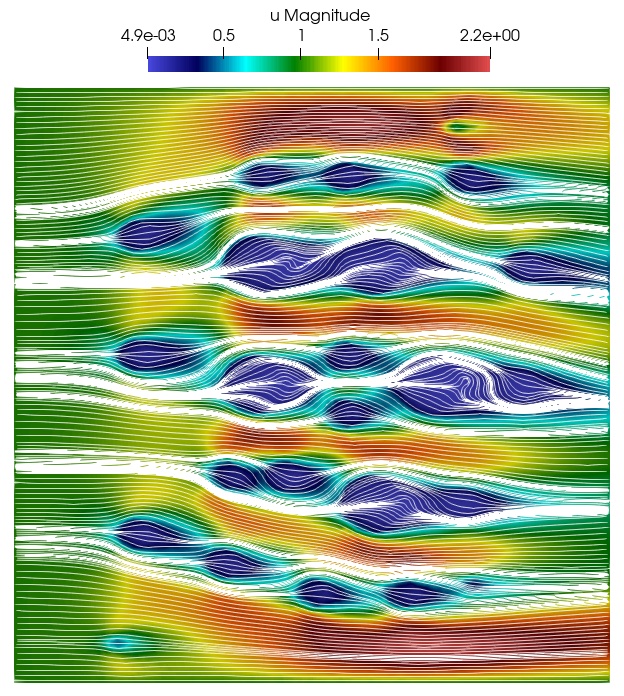}\\
\includegraphics[width=0.78\linewidth]{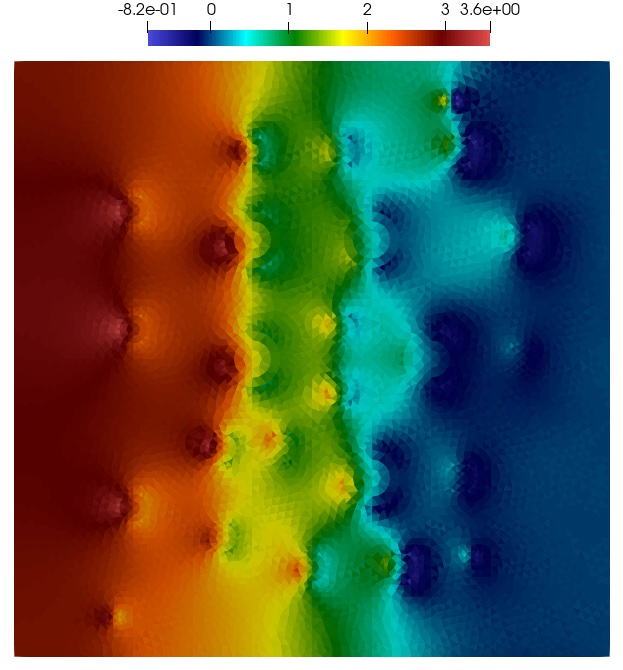}
\end{subfigure}
\begin{subfigure}{0.32\textwidth}
\centering
\includegraphics[width=0.78\linewidth]{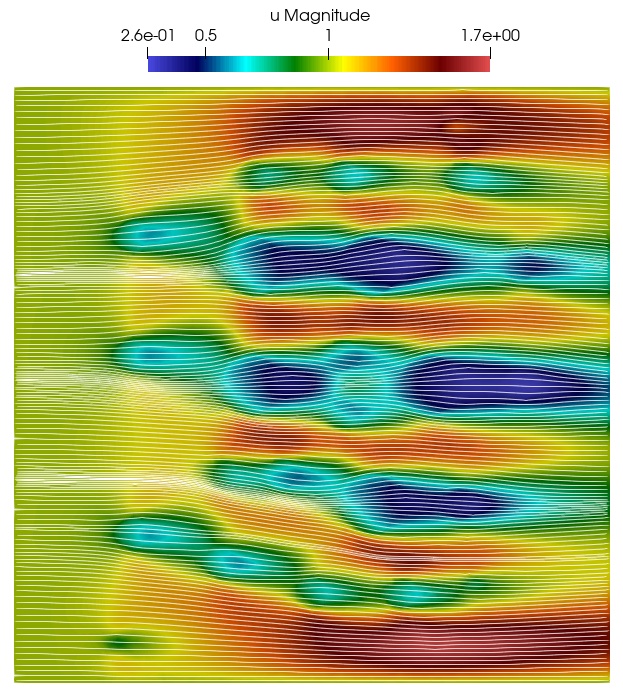}\\
\includegraphics[width=0.78\linewidth]{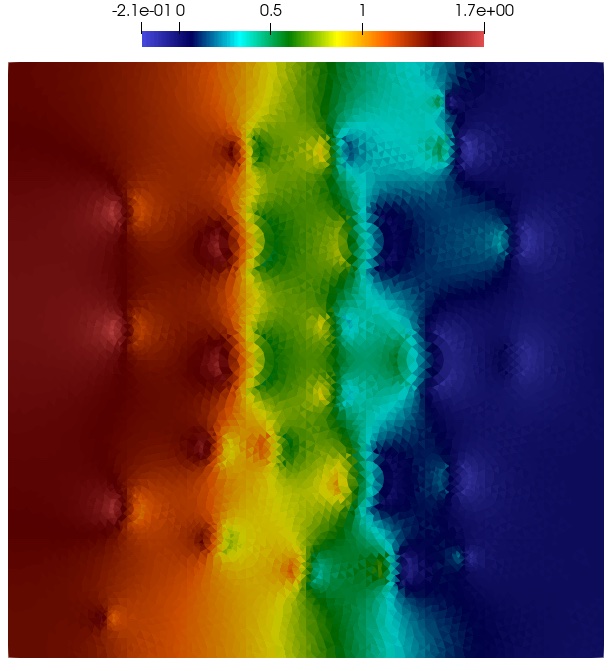}
\end{subfigure}
\caption{\textit{Test 3} ($Re = 100$, $C = 1$).  Fine grid solution at final time for $Da = 10^{-5}, 10^{-4}, 10^{-3}$ (from left to right).  
First row: velocity magnitude with streamlines. 
Second row: pressure}
\label{up-3}
\end{figure}

The fine grid numerical solutions computed at the final time step for Tests 1, 2, and 3 are displayed in Figures \ref{up-1}, \ref{up-2}, and \ref{up-3}, respectively. Our numerical implementation has been performed using FEniCS finite element library~\cite{logg2012automated}. Paraview software \cite{ahrens2005paraview} has been utilized for the visualisation of the results. The velocity field magnitudes along with 100 streamlines are exhibited in the first row (of each figure) for three different Darcy numbers, $Da = 10^{-5}, 10^{-4}$ and $10^{-3}$. The computed pressure fields corresponding to these Darcy number values are shown in the respective second rows. Note that the smallest $Da = 10^{-5}$ corresponds to almost impermeable inclusions case while $Da = 10^{-3}$ represents the situation with highly permeable inclusions. We observe that the Darcy number has a significant impact on the two-dimensional flow velocity and pressure fields with multiple inclusions in all three Test cases. It may be noted that larger Darcy number (less permeable case) yields smaller velocity magnitudes and pressures.

Typical flow field scenarios emerging from the fine grid solutions, computed in our three numerical Test runs, may be described as follows. As seen in Figure \ref{up-1} for Test 1, the effect of the subdomains $\Omega_f$, that is the porous inclusions, is to partially divert the fluid flow. Many streamlines bypass the inclusions, whereas other pass through small parts. When the Darcy number is small, the velocity streamline penetration from $\Omega_f$ is less for small Reynolds number and Forchheimer coefficient values ($Re =1$ and $C=1$, Figure \ref{up-1} top row left). But the instantaneous fluid streamlines inside $\Omega_p$ appear to undergo significant changes. In particular, the streamlines bend towards the center of inclusions as seen in the figure. This is probably due to the reason that the boundary effects are stronger in the case of nearly impermeable inclusions. The fluid penetration gets higher as the Darcy number increases as in Figure \ref{up-1} (top row middle and right). The instantaneous streamlines from $\Omega_f$ penetrate into $\Omega_f$ relatively easier when the Darcy number is high. The fluid lines in the interior domain follows approximate straight line paths in this case. The boundary effects become weaker when the Darcy number gets higher. Indeed, the presence of more void volume in the inclusions cause more fluid flow thorough $\Omega_p$. Also, the effect of nonlinearity is weaker since the porous inclusions allow a finite amount of fluid pass through with non-zero velocity at the interfaces. The pressure increases from left to right as shown in Figure \ref{up-1} bottom row. Relative changes in the high pressure zones in the exterior domain $\Omega_p$ can be noticed as the Darcy number increases.       

The streamlines (along with velocity magnitudes) and the pressure fields for Test 2 with $Re =10$ and $C=10$ are presented in Figure \ref{up-2} for the indicated Darcy numbers. The changes in the pattern is generally similar to that in Test 1 case. But in the present case, the fluid velocity appears to be vigorous in porous domains $\Omega_p$. Bending of streamline patterns  appear to increase for $Da = 10^{-5}$ (low Darcy number) as can be noticed in this figure. This could possibly be due to the effect of high values of Reynolds and Forchheimer numbers. The nonlinear convective terms in both $\Omega_p$ and $\Omega_f$ impact the velocity changes inside the porous domains. The numerical values for the pressure shows a sort of decreasing trend (Figure \ref{up-2} bottom row) in comparison with Test 1 case results.   

Numerical simulation results for Test 3 with $Re =100$ and $C=1$ are portrayed in Figure \ref{up-3} for three Darcy numbers. It can be seen that  for a relatively high Reynolds number with low $Da$ influence the flow structures in both $\Omega_f$ and $\Omega_p$ significantly. The recirculating zones/wakes in both $\Omega_f$ and $\Omega_f$ domains are seen and secondary flow patterns develop at the rear exit of the porous inclusions. This scenario may be expected since the Reynolds number in the exterior domain $\Omega_f$ is much higher than in Tests 1 and 2. The effect of nonlinearity can be recognised from these flow topologies generated via fine grid numerical solutions. The pressure decreases as the Darcy number increases as noticed from Figure \ref{up-3} (bottom row).   

\begin{table}[h!]
\center
\begin{tabular}{|c|c|ccc|ccc|}
\hline
 & & 
 \multicolumn{3}{|c|}{without oversampling} &
  \multicolumn{3}{|c|}{with oversampling}\\ 
 \raisebox{1.5ex}[0cm][0cm]{$\mathcal{M}$ } 
&\raisebox{1.5ex}[0cm][0cm]{$DOF_H$ } 
& $e_u$  & $e_s$ & $e_p$ 
& $e_u$  & $e_s$ & $e_p$ \\
\hline
\multicolumn{8}{|c|}{$Da = 10^{-5}$} \\
\hline
5 & 600 	& 29.291 	& 82.109 	& $>$100 	& 50.039 	& 97.617 	& $>$100 \\
10 & 1100 	& 20.826 	& 63.059 	& $>$100 	& 32.391 	& 66.400 	& $>$100 \\
15 & 1600 	& 11.894 	& 44.644 	& 45.406 	& 6.781 	& 19.869 	& 6.083 \\
20 & 2100 	& 7.243 	& 35.325 	& 32.736 	& 3.150 	& 10.224 	& 0.288 \\
25 & 2600 	& 4.633 	& 29.107 	& 25.135 	& 1.300 	& 3.940 	& 1.409 \\
\hline
\multicolumn{8}{|c|}{$Da = 10^{-4}$} \\
\hline
5 & 600 	& 32.946 	& 88.932 	& $>$100 	& 30.780 	& 84.187 	& $>$100 \\
10 & 1100 	& 19.946 	& 64.455 	& $>$100 	& 11.542 	& 36.650 	& 34.310 \\
15 & 1600 	& 12.247 	& 49.235 	& 46.594 	& 3.676 	& 14.492 	& 2.427 \\
20 & 2100 	& 6.412 	& 33.629 	& 24.532 	& 1.616 	& 6.813 	& 0.547 \\
25 & 2600 	& 5.346 	& 31.003 	& 21.447 	& 0.921 	& 3.432 	& 1.063 \\
\hline
\multicolumn{8}{|c|}{$Da = 10^{-3}$} \\
\hline
5 & 600 	& 19.650 	& 90.158 	& 94.043 	& 18.440 	& 86.225 	& $>$100 \\
10 & 1100 	& 11.072 	& 70.370 	& 30.690 	& 2.369 	& 19.114 	& 2.778 \\
15 & 1600 	& 6.543 	& 53.626 	& 18.934 	& 1.113 	& 9.963 	& 0.588 \\
20 & 2100 	& 3.658 	& 37.045 	& 12.094 	& 0.469 	& 4.400 	& 0.268 \\
25 & 2600 	& 3.130 	& 33.896 	& 10.657 	& 0.342 	& 3.184 	& 0.384 \\
\hline
\end{tabular}
\caption{\textit{Test 1} ($Re = 1$, $C = 1$). Relative errors for velocity, stress and pressure.   
}
\label{table-1}
\end{table}

\begin{table}[h!]
\center
\begin{tabular}{|c|c|ccc|ccc|}
\hline
 & & 
 \multicolumn{3}{|c|}{without oversampling} &
  \multicolumn{3}{|c|}{with oversampling}\\ 
 \raisebox{1.5ex}[0cm][0cm]{$\mathcal{M}$ } 
&\raisebox{1.5ex}[0cm][0cm]{$DOF_H$ } 
& $e_u$  & $e_s$ & $e_p$ 
& $e_u$  & $e_s$ & $e_p$ \\
\hline
\multicolumn{8}{|c|}{$Da = 10^{-5}$} \\
\hline
5 & 600 	& 29.277 	& 81.882 	& $>$100 	& 51.075 	& 97.439 	& $>$100 \\
10 & 1100 	& 21.291 	& 63.432 	& $>$100 	& 22.844 	& 55.525 	& 73.378 \\
15 & 1600 	& 12.149 	& 44.955 	& 45.231 	& 6.805 	& 18.657 	& 5.607 \\
20 & 2100 	& 7.419 	& 35.625 	& 32.596 	& 2.628 	& 8.167 	& 0.517 \\
25 & 2600 	& 4.837 	& 29.392 	& 25.121 	& 1.370 	& 3.968 	& 1.434 \\
\hline
\multicolumn{8}{|c|}{$Da = 10^{-4}$} \\
\hline
5 & 600 	& 31.062 	& 86.300 	& $>$100 	& 30.716 	& 83.064 	& $>$100 \\
10 & 1100 	& 20.410 	& 64.989 	& $>$100 	& 13.897 	& 40.996 	& 47.419 \\
15 & 1600 	& 12.828 	& 49.372 	& 51.047 	& 4.220 	& 15.734 	& 3.668 \\
20 & 2100 	& 6.667 	& 33.709 	& 26.018 	& 1.521 	& 6.280 	& 0.766 \\
25 & 2600 	& 5.558 	& 30.790 	& 22.550 	& 1.038 	& 3.609 	& 1.161 \\
\hline
\multicolumn{8}{|c|}{$Da = 10^{-3}$} \\
\hline
5 & 600 	& 24.552 	& 87.678 	& $>$100 	& 23.846 	& 86.250 	& $>$100 \\
10 & 1100 	& 15.524 	& 68.307 	& 78.156 	& 4.944 	& 24.293 	& 9.542 \\
15 & 1600 	& 8.945 	& 51.202 	& 36.681 	& 2.503 	& 13.725 	& 1.662 \\
20 & 2100 	& 4.783 	& 35.877 	& 20.730 	& 0.923 	& 5.508 	& 0.562 \\
25 & 2600 	& 4.087 	& 33.537 	& 18.649 	& 0.607 	& 3.460 	& 0.784 \\
\hline
\end{tabular}
\caption{\textit{Test 2} ($Re = 10$, $C = 10$). Relative errors for velocity, stress and pressure.   
}
\label{table-2}
\end{table}

\begin{table}[h!]
\center
\begin{tabular}{|c|c|ccc|ccc|}
\hline
 & & 
 \multicolumn{3}{|c|}{without oversampling} &
  \multicolumn{3}{|c|}{with oversampling}\\ 
 \raisebox{1.5ex}[0cm][0cm]{$\mathcal{M}$ } 
&\raisebox{1.5ex}[0cm][0cm]{$DOF_H$ } 
& $e_u$  & $e_s$ & $e_p$ 
& $e_u$  & $e_s$ & $e_p$ \\
\hline
\multicolumn{8}{|c|}{$Da = 10^{-5}$} \\
\hline
5 & 600 	& 40.209 	& 89.317 	& $>$100 	& 55.351 	& 97.739 	& $>$100 \\
10 & 1100 	& 29.889 	& 71.256 	& 79.574 	& 32.087 	& 69.818 	& $>$100 \\
15 & 1600 	& 16.864 	& 53.072 	& 29.122 	& 12.628 	& 31.352 	& 4.277 \\
20 & 2100 	& 11.126 	& 44.202 	& 21.655 	& 5.819 	& 17.392 	& 5.428 \\
25 & 2600 	& 8.491 	& 37.592 	& 15.085 	& 5.085 	& 15.746 	& 5.644 \\
\hline
\multicolumn{8}{|c|}{$Da = 10^{-4}$} \\
\hline
5 & 600 	& 33.282 	& 88.705 	& $>$100 	& 37.595 	& 88.684 	& $>$100 \\
10 & 1100 	& 24.283 	& 69.830 	& 88.757 	& 20.383 	& 55.240 	& 64.051 \\
15 & 1600 	& 14.044 	& 52.181 	& 34.147 	& 9.611 	& 29.225 	& 3.453 \\
20 & 2100 	& 8.653 	& 39.735 	& 17.283 	& 4.112 	& 16.339 	& 2.690 \\
25 & 2600 	& 7.646 	& 36.807 	& 15.804 	& 3.383 	& 14.139 	& 3.024 \\
\hline
\multicolumn{8}{|c|}{$Da = 10^{-3}$} \\
\hline
5 & 600 	& 20.627 	& 93.689 	& $>$100 	& 20.836 	& 90.268 	& $>$100 \\
10 & 1100 	& 14.624 	& 74.773 	& 89.399 	& 8.856 	& 45.058 	& 24.358 \\
15 & 1600 	& 7.288 	& 50.458 	& 34.458 	& 5.089 	& 31.344 	& 5.382 \\
20 & 2100 	& 5.192 	& 42.230 	& 20.437 	& 1.511 	& 13.154 	& 0.477 \\
25 & 2600 	& 4.660 	& 39.015 	& 18.250 	& 1.230 	& 10.482 	& 1.095 \\
\hline
\end{tabular}
\caption{\textit{Test 3} ($Re = 100$, $C = 1$).  Relative errors for velocity, stress and pressure.  
}
\label{table-3}
\end{table}

\begin{figure}[h!]
\centering
\begin{subfigure}{0.32\textwidth}
\centering
Reference solution
\end{subfigure}
\begin{subfigure}{0.32\textwidth}
\centering
MS without oversampling\\
$e_u = 4.6 \%$,  $e_p = 25.1 \%$
\end{subfigure}
\begin{subfigure}{0.32\textwidth}
\centering
MS with oversampling\\
$e_u = 1.3 \%$,  $e_p = 1.4 \%$
\end{subfigure}\\
\begin{subfigure}{0.32\textwidth}
\centering
\includegraphics[width=0.78\linewidth]{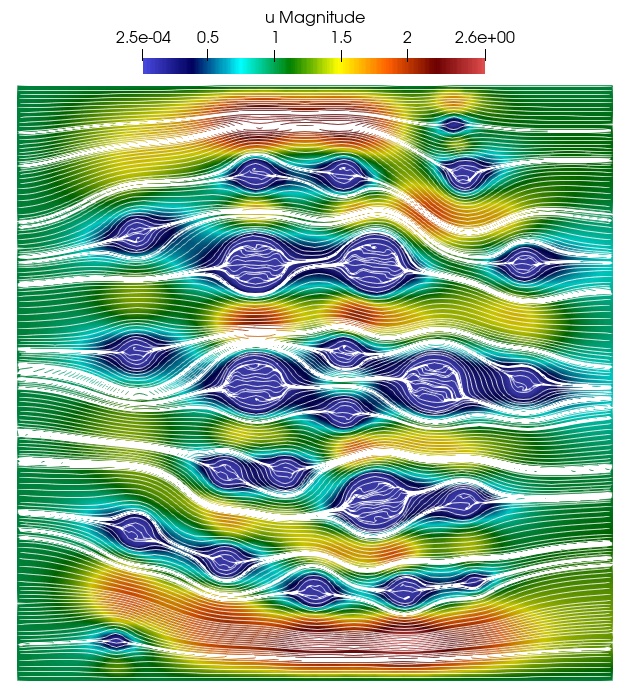}\\
\includegraphics[width=0.78\linewidth]{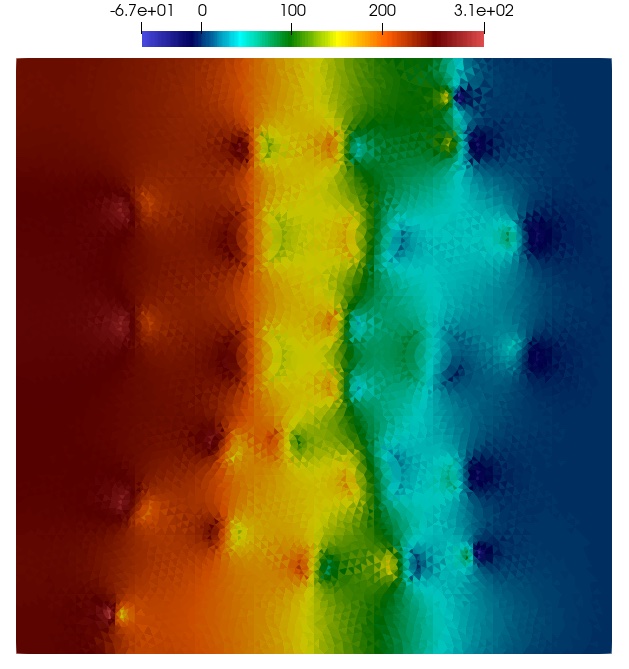}
\end{subfigure}
\begin{subfigure}{0.32\textwidth}
\centering
\includegraphics[width=0.78\linewidth]{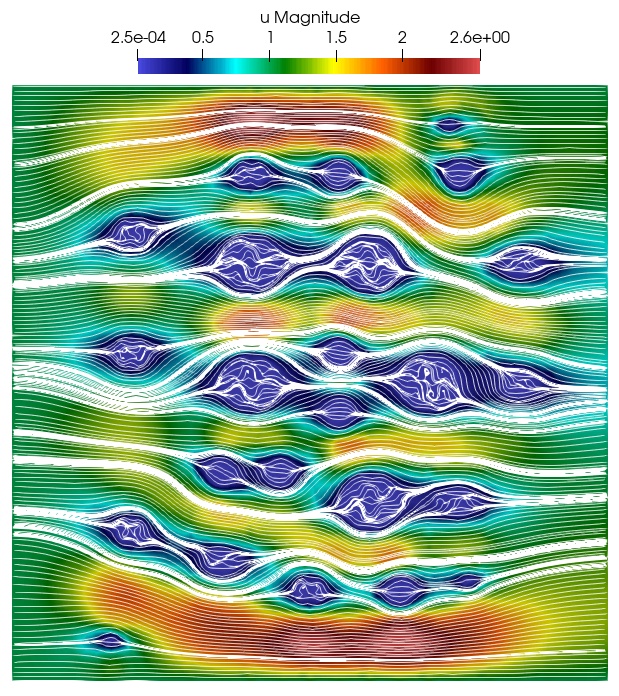}\\
\includegraphics[width=0.78\linewidth]{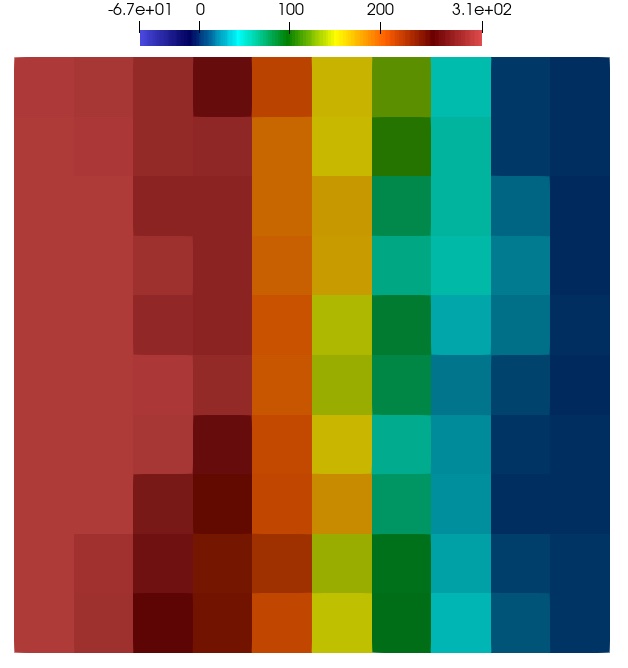}
\end{subfigure}
\begin{subfigure}{0.32\textwidth}
\centering
\includegraphics[width=0.78\linewidth]{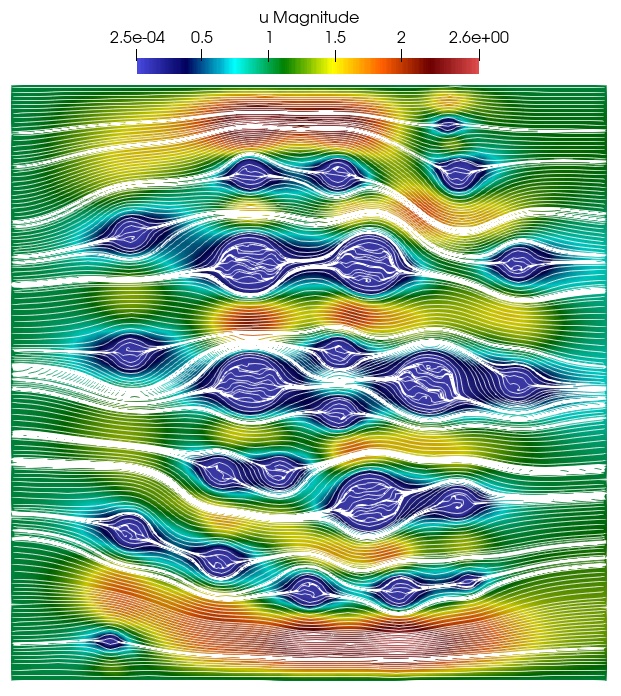}\\
\includegraphics[width=0.78\linewidth]{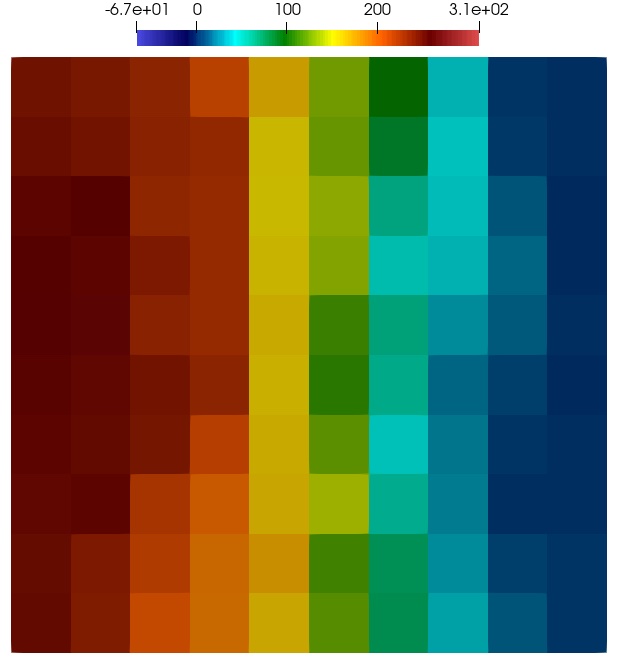}
\end{subfigure}
\caption{\textit{Test 1} ($Re = 1$, $C = 1$,  $Da = 10^{-5}$).   
First row: velocity magnitude with streamlines. 
Second row: pressure. 
Reference solution: $DOF_h = 87528$.  
Multiscale solution: $DOF_H=2600$}
\label{up-ms-1}
\end{figure}

\begin{figure}[h!]
\centering
\begin{subfigure}{0.32\textwidth}
\centering
Reference solution
\end{subfigure}
\begin{subfigure}{0.32\textwidth}
\centering
MS without oversampling\\
$e_u = 4.8 \%$,  $e_p = 25.1 \%$
\end{subfigure}
\begin{subfigure}{0.32\textwidth}
\centering
MS with oversampling\\
$e_u = 1.3 \%$,  $e_p = 1.4 \%$
\end{subfigure}\\
\begin{subfigure}{0.32\textwidth}
\centering
\includegraphics[width=0.78\linewidth]{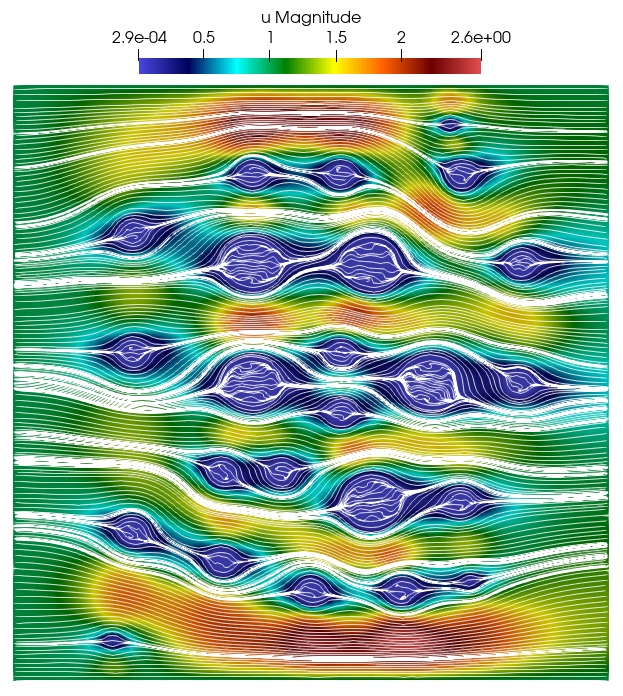}\\
\includegraphics[width=0.78\linewidth]{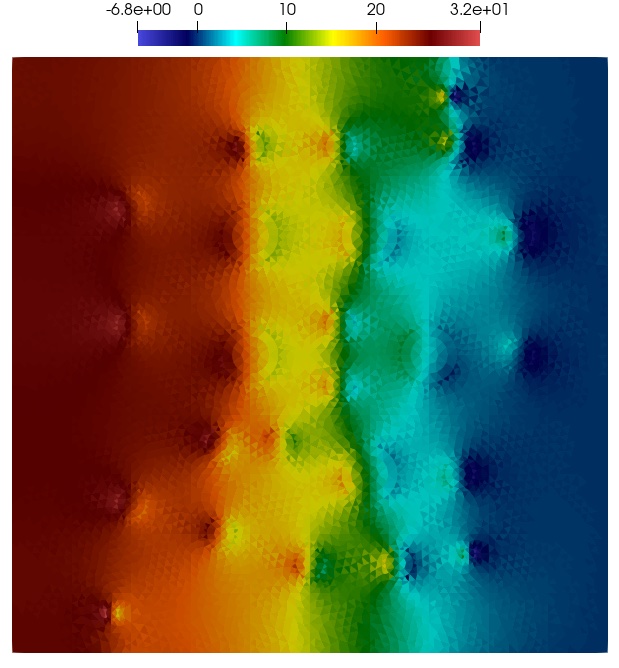}
\end{subfigure}
\begin{subfigure}{0.32\textwidth}
\centering
\includegraphics[width=0.78\linewidth]{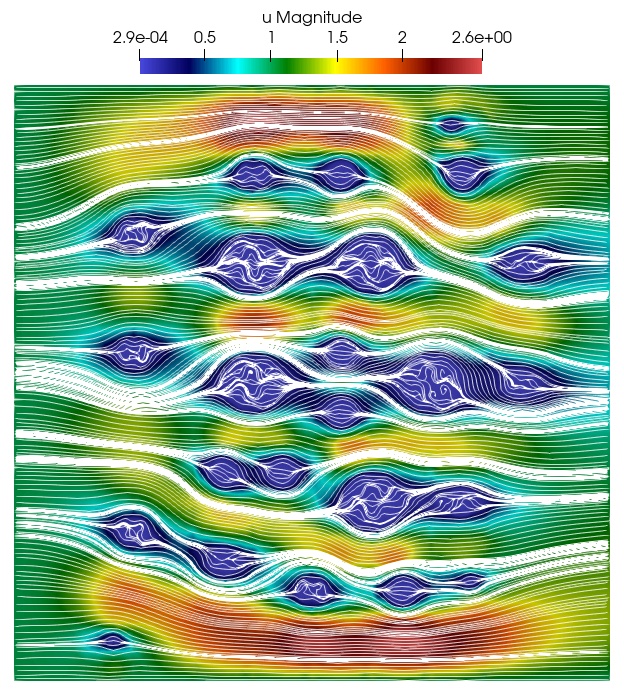}\\
\includegraphics[width=0.78\linewidth]{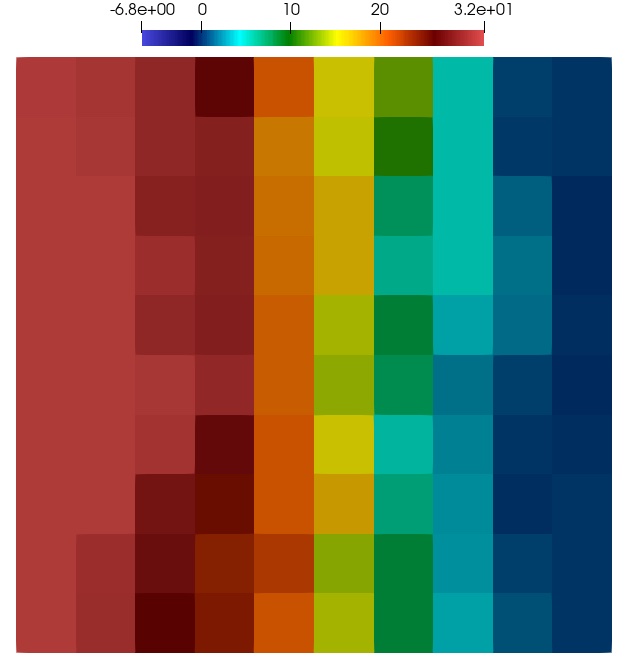}
\end{subfigure}
\begin{subfigure}{0.32\textwidth}
\centering
\includegraphics[width=0.78\linewidth]{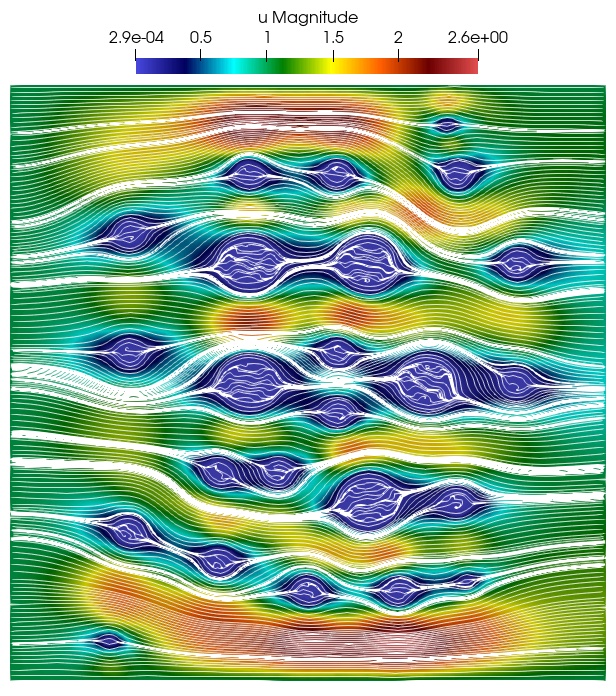}\\
\includegraphics[width=0.78\linewidth]{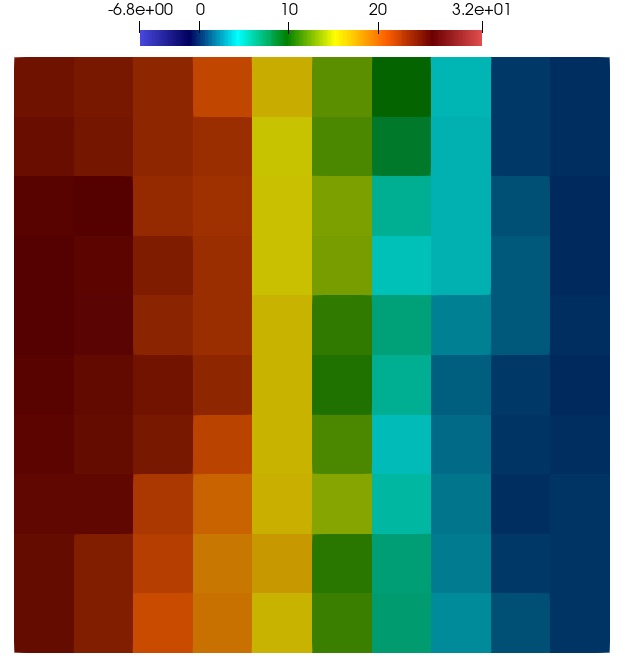}
\end{subfigure}
\caption{\textit{Test 2} ($Re = 10$, $C = 10$,  $Da = 10^{-5}$).   
First row: velocity magnitude with streamlines. 
Second row: pressure.  
Reference solution: $DOF_h = 87528$.  
Multiscale solution: $DOF_H=2600$}
\label{up-ms-2}
\end{figure}

\begin{figure}[h!]
\centering
\begin{subfigure}{0.32\textwidth}
\centering
Reference solution
\end{subfigure}
\begin{subfigure}{0.32\textwidth}
\centering
MS without oversampling\\
$e_u = 8.4 \%$,  $e_p = 15.0 \%$
\end{subfigure}
\begin{subfigure}{0.32\textwidth}
\centering
MS with oversampling\\
$e_u = 5.0 \%$,  $e_p = 5.6 \%$
\end{subfigure}\\
\begin{subfigure}{0.32\textwidth}
\centering
\includegraphics[width=0.78\linewidth]{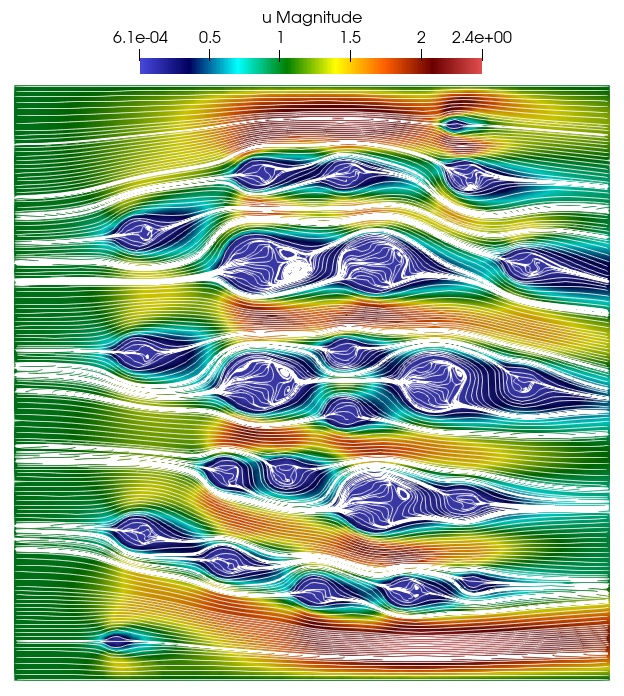}\\
\includegraphics[width=0.78\linewidth]{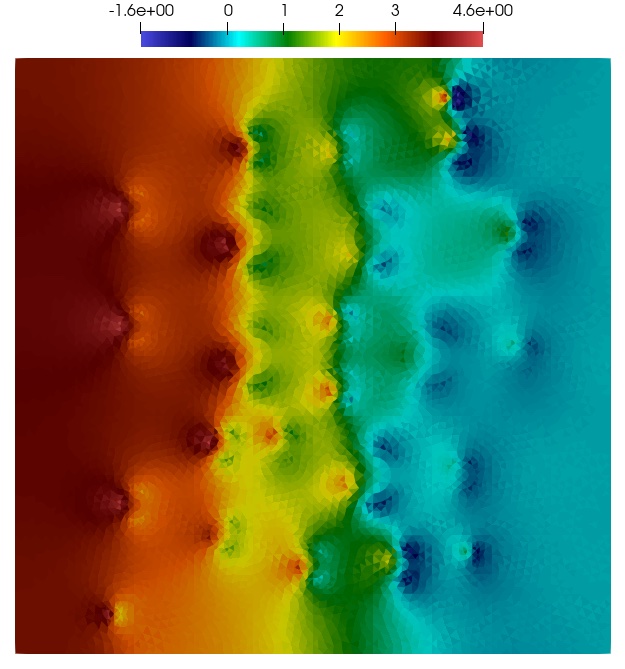}
\end{subfigure}
\begin{subfigure}{0.32\textwidth}
\centering
\includegraphics[width=0.78\linewidth]{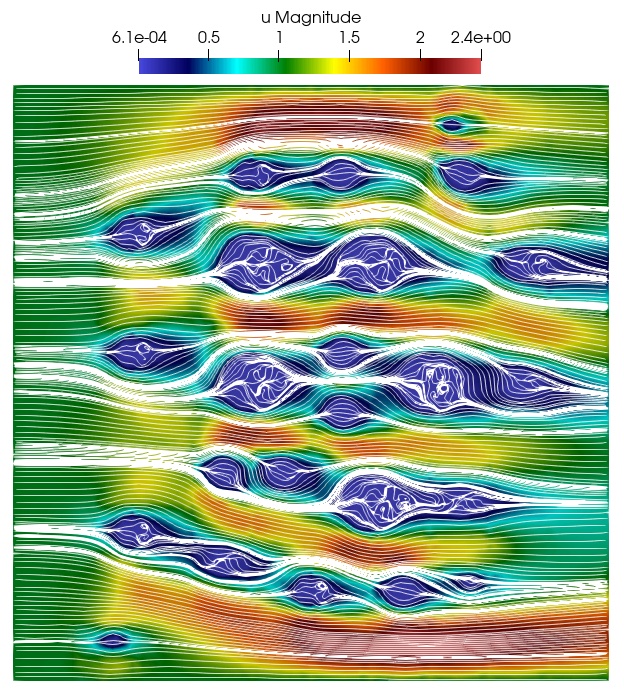}\\
\includegraphics[width=0.78\linewidth]{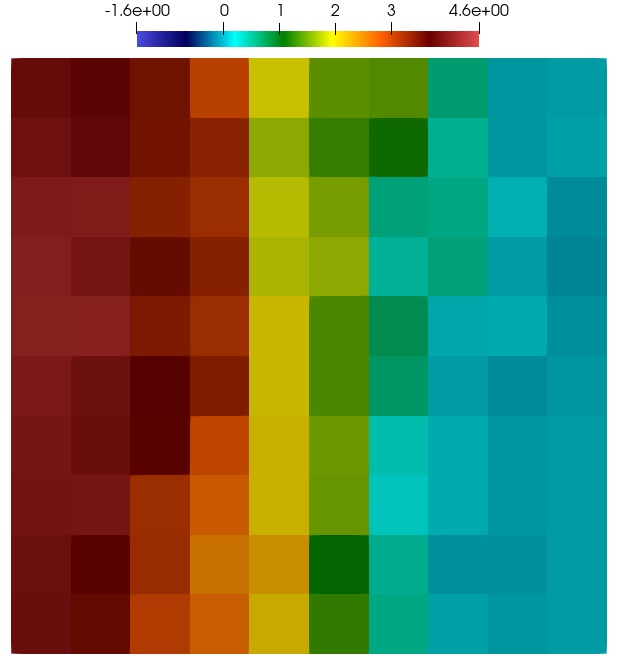}
\end{subfigure}
\begin{subfigure}{0.32\textwidth}
\centering
\includegraphics[width=0.78\linewidth]{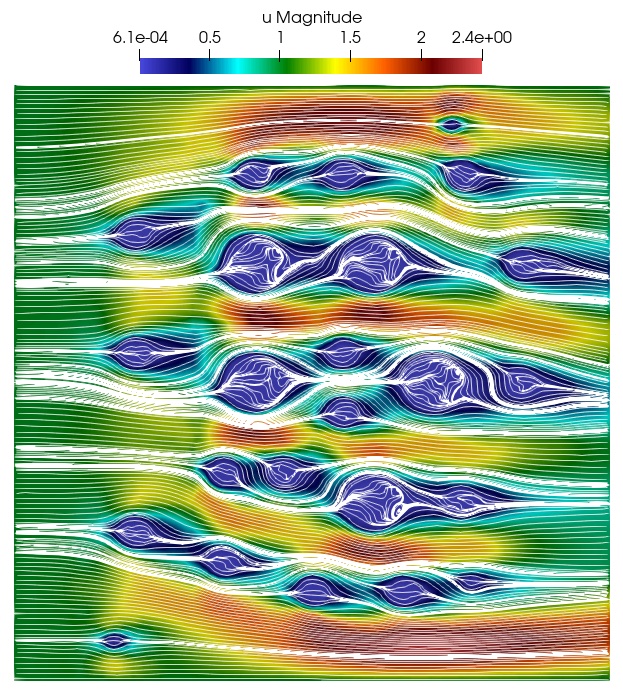}\\
\includegraphics[width=0.78\linewidth]{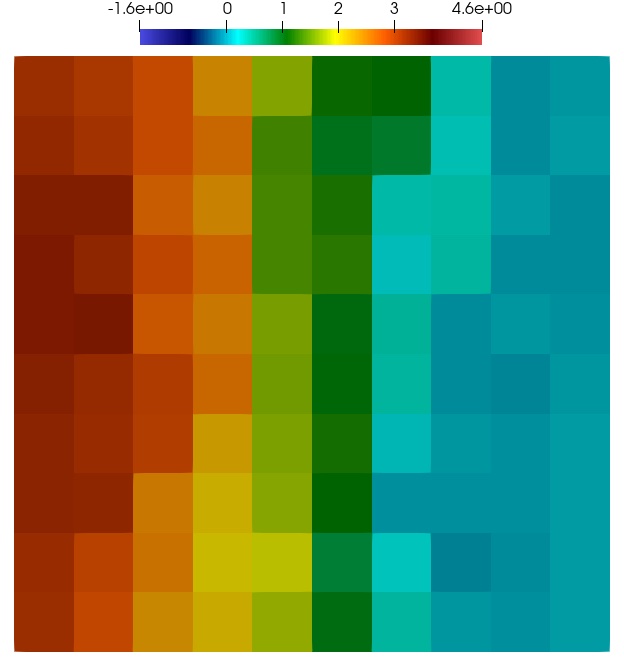}
\end{subfigure}
\caption{\textit{Test 3} ($Re = 100$, $C = 1$,  $Da = 10^{-5}$).   
First row: velocity magnitude with streamlines. 
Second row: pressure.  
Reference solution: $DOF_h = 87528$.  
Multiscale solution: $DOF_H=2600$}
\label{up-ms-3}
\end{figure}

Next, we discuss numerical solutions of the same problem computed using multisolver (MS). In this study we perform numerical simulations by varying the number of multiscale basis functions $\mathcal{M}$ ranging from 5 to 25. Additionally, we compute results with and without oversampling strategy in basis constructions. Also, for comparison purposes we use the fine grid solution as a reference solution to the problem. Note that the size of the system on fine grid is chosen to be $DOF_h = 87528$.

To compare the results, we calculate the relative $L_2$ error between reference solution $(u_h, p_h)$ and multiscale solutions $(u_{ms}, p_{ms})$  for stress ($e_s$), velocity ($e_u$) and pressure ($e_p$) using  
\[
e_s = \left( \frac{ 
\int_{\Omega} \epsilon(u_h - u_{ms}) \colon \epsilon(u_h - u_{ms}) dx }{
\int_{\Omega} \epsilon(u_h) \colon \epsilon(u) dx} \right)^{1/2}, 
\]\[
e_u = \left( \frac{ 
\int_{\Omega} (u_h - u_{ms}) \cdot (u - u_{ms}) \ dx }{
\int_{\Omega} u_h \cdot u_h \ dx } \right)^{1/2}, 
\quad 
e_p = \left( \frac{ 
\int_{\Omega} (\overline{p}_h - p_{ms})^2 \ dx }{
\int_{\Omega} \overline{p}_h^2 \ dx }  \right)^{1/2}, 
\]
where $\epsilon(u) = (\nabla u + \nabla u^T)/2$ is the strain tensor. We note that for pressure, we use $L_2$ errors on the coarse grid and $\overline{p}$ denotes the coarse cell average for reference (fine-grid) pressure. 

The computed relative errors for stress, velocity and pressure at the final time for Tests 1, 2 and 3 are provided in Tables \ref{table-1}, \ref{table-2} and \ref{table-3}, respectively. The results are given for three different values of the Darcy number $Da = 10^{-5}, 10^{-4}$ and $10^{-3}$. The cases with and without oversampling strategy in basis construction have been considered in solving the equations using multiscale solver.  In order to construct oversampled region, we  used one additional coarse grid layer in local domain construction (see Figure \ref{meshc}). The first column in these tables shows the number of multiscale basis functions $\mathcal{M}$ for the velocity in each local domain, the second column specifies the dimension of the multiscale space ($DOF_H$), and the remaining columns indicate the relative errors in percentage. It can be seen from the tables that reasonably good solutions are obtained using sufficient number of $\mathcal{M}$ for accurate approximation of the velocity field. Note that the use of 5 or 10 multiscale basis functions in our multiscale approach is not sufficient for computing reasonable numerical solutions for both test problems. 
For instance, from Table \ref{table-1} for Test 1 with $Da = 10^{-5}$, we obtain $20.8$ and $32.3$ \% of relative velocity errors using 10 multiscale basis functions for the case without and with oversampling strategy.  However, when we use 25 multiscale basis functions, we have $4.6$ and $1.3$ \% of relative  errors for velocity field. This implies that significant relative error reductions for the physical quantities can be achieved by employing additional multiscale basis functions. The higher impact of the oversampling strategy on the multiscale method errors can also be observed in the table of values.  Using oversampled local domains in multiscale basis space construction, one can obtain velocity and pressure fields with the margin of errors around $1$ \% by using a sufficiently large number of multiscale basis functions.

In Tables \ref{table-2} and  \ref{table-3}, we provide results for the relative errors in velocity, stress and pressure for the test problems with $Re = C =10$ (Test 2) and $Re = 100$, $C = 1$ (Test 3). Note that the nonlinear effects are significant since the Reynolds number is higher in these cases. Observe that the errors are larger than those found for Test 1 with $Re = C = 1$.  However, we still obtain decent numerical results with $1$\% of velocity error using the multiscale method with oversampling strategy for Test 2 and $1\%-5\%$ of velocity error for Test 3, respectively. 
We also notice that the oversampling strategy has a huge impact on the pressure accuracy as well. For example, in Test 2 with $Da = 10^{-5}$, we have $4.8$ and $25.1$ \% of velocity and pressure errors without oversampling strategy for $\mathcal{M} = 25$.  By applying the oversampling approach, the velocity and pressure relative errors reduce to $1.3$ and $1.4$ \%.  
The influence of Darcy number on the multiscale method accuracy may also be noticed from Tables \ref{table-1},  \ref{table-2} and \ref{table-3}.  The relative error is smaller for high Darcy number flows (fluid flows with highly permeable inclusions) in all test problems, but it is larger for less permeable inclusions (low Darcy number flows). In particular, we have  $1.3$ \% of errors for velocity  and $1.4$ \% for pressure in Test 1 with $Da = 10^{-5}$ using 25 multiscale basis functions with oversampling. For $Da=10^{-3}$, we obtain relative errors $0.3$ \% for velocity  and $0.3$ \% for pressure with $\mathcal{M}=25$.  
For the Test 3, we obtain $5.0$ \% of velocity error for $Da=10^{-5}$ and $1.2$ \% for $Da=10^{-3}$.  Note that the permeability of the inclusions, alternatively the $Da$, has a greater impact on the  velocity field, where for less permeable case we obtain more heterogeneous velocity field  that is harder to approximate. This indicated that more multiscale velocity basis functions may be needed for the mutiscale solver to capture such velocity fields associated with complex flows.  
We also observe the significant influence of the nonlinear parameters ($Re$ and $C$) linked to the flow in heterogeneous media (see Table \ref{table-3}).  Our multiscale solver results yield larger errors in the case of very high Reynolds numbers. This is perhaps due to the way of our multiscale basis functions construction in the adopted multiscale approach. We remark that in the present study we have utilised the linear basis construction based on the global velocity field. One can apply the online approach that takes the residual into account while constructing the basis functions \cite{chung2016online} to achieve significant reduction of errors, especially in high Reynolds number flow situations. But we will not address those techniques in this paper.

For comparison, the flow fields and the pressure obtained based on the fine grid (reference) solution and multiscale solution with 25 basis functions are shown in Figures \ref{up-ms-1}, \ref{up-ms-2} and \ref{up-ms-3} for Test 1, 2 and 3 with $Da = 10^{-5}$. 
The velocity magnitude along with streamlines and the pressure field for the reference solution on the fine grid  are displayed in the first column of the figures. The respective physical quantities based on multiscale solutions without oversampling strategy are shown in second column and results using oversampling  strategy are displayed in the third column. The number of multiscale functions is chosen to be $\mathcal{M}=25$.   We observe that the oversampling approach provides a large error reduction in the physical quantities.  Moreover, by comparing results from second and third columns for the velocity field, we notice that for the case without oversampling approach the errors are concentrated mostly near local domain boundary.  By applying oversampling strategy we can reduce the boundary effects in snapshot space construction and  improve the accuracy of multiscale methods. A close agreement of the flow features between fine grid and multiscale solutions is apparent in these figures.

\section{Conclusion}

A powerful multiscale solver (multisolver)  for computing solutions of the Navier-Stokes/Darcy-Brinkman-Forchheimer (NSDBM) model problem describing two-dimensional flows with several circular porous inclusions is presented in this work. A scheme for the construction of multiscale velocity basis functions in heterogeneous domains with permeable obstacles by the use of GMsFEM framework is narrated. Specifically, in order to construct the basis functions, one generates a snapshot space in local domains with and without oversampling strategy. The snapsots are actually obtained via the solutions of the local problems with all possible boundary conditions. One then solves a local spectral problem to reduce the size of the snapshot space and constructs a low dimensional local representation of the solution. Numerical results are reported for three different test cases for the choices chosen for Reynolds number and Forchheimer coefficient to demonstrate the power of our method. Multiscale solutions are also displayed for various Darcy numbers to show the influence of inclusions permeability and the multiscale performance. The following points are noteworthy from the present numerical study.
\begin{itemize}
	\item Increase in the number of multiscale basis functions ($\mathcal{M}$) and oversampling technique enhances the accuracy of the velocity magnitudes and pressure. 
	\item The discrete system size for the mutiscale method used is approximately 34 times smaller than that for fine grid set up. Thus, there is a significant dimension reduction in the multiphase problem under consideration, elucidating the power of multiscale approach.
	\item The choice $\mathcal{M} = 25$ yields numerical results with relative errors around $1\%$. The method indicates that larger errors occur close to the boundary of the inclusions when $Re = 100$. Note that in this case the nonlinearity contribution is rather significant.
	\item For small Reynolds number flows with $Re = 1, 10$ with oversampling strategy, the performance of the proposed method is very good.
	\item Our numerical solutions capture the flow features with reasonable accuracy. Depicted instantaneous streamline topologies show recirculatory flow patterns and attached wakes at the rear exit of the porous inclusions for low permeability (low Darcy number). Stronger fluid penetration from $\Omega_f$ into permeable subdomains $\Omega_p$ is observed for higher Darcy numbers (high permeable inclusions).   
	   
\end{itemize}

It should be pointed out that we have considered two dozen heterogeneous porous inclusions in this investigation mainly for testing our numerical approach and to illustrate typical flow pattern scenarios. However, the method can handle flow problems with any number of permeable inclusions with arbitrary locations. Extensive study of situations by varying the key parameters in physically admissible ranges may reveal flow transitions and further development of recirculating zones and wakes in the clear fluid domain $\Omega_f$ as well as in porous inclusions $\Omega_p$. Our multisolver solutions and simulations presented herein may be crucial in such exhaustive studies. Finally, it may be worthwhile to mention that one can apply online approach discussed in \cite{chung2015residual, chung2016online} to achieve remarkable error reduction in the numerical solutions for the physical quantities. 

\bibliographystyle{plain}
\bibliography{lit3}

\end{document}